\input  psfig.sty
%
%
\catcode`@=11 
%
%

\def\b@lank{ }

\newif\if@simboli
\newif\if@riferimenti
\newif\if@bozze

\def\bozze{\@bozzetrue\font\tt@bozze=cmtt8}
\def\og@gi{\number\day\space\ifcase\month\or 
   gennaio\or febbraio\or marzo\or aprile\or maggio\or giugno\or 
   luglio\or agosto\or settembre\or ottobre\or novembre\or dicembre\fi
   \space\number\year}
\newcount\min@uti
\newcount\or@a
\newcount\ausil@iario
\min@uti=\number\time
\or@a=\number\time
\divide\or@a by 60
\ausil@iario=-\number\or@a
\multiply\ausil@iario by 60
\advance\min@uti by \number\ausil@iario
\def\ora@esecuzione{\the\or@a:\the\min@uti}  
\def\makefootline{\baselineskip=24pt\line{\the\footline}
    \if@bozze\vskip-10pt\tt@bozze
             \noindent \jobname\hfill\og@gi, ore \ora@esecuzione\fi}

\newwrite\file@simboli
\def\simboli{
    \immediate\write16{ !!! Genera il file \jobname.SMB }
    \@simbolitrue\immediate\openout\file@simboli=\jobname.smb}

\newwrite\file@ausiliario
\def\riferimentifuturi{
    \immediate\write16{ !!! Genera il file \jobname.AUX }
    \@riferimentitrue\openin1 \jobname.aux
    \ifeof1\relax\else\closein1\relax\input\jobname.aux\fi
    \immediate\openout\file@ausiliario=\jobname.aux}

\newcount\eq@num\global\eq@num=0
\newcount\sect@num\global\sect@num=0

\newif\if@ndoppia
\def\numerazionedoppia{\@ndoppiatrue\gdef\la@sezionecorrente{\the\sect@num}}

\def\se@indefinito#1{\expandafter\ifx\csname#1\endcsname\relax}
\def\spo@glia#1>{} 

\newif\if@primasezione
\@primasezionetrue

\def\s@ection#1\par{\immediate
    \write16{#1}\if@primasezione\global\@primasezionefalse\else\goodbreak
    \vskip\spaziosoprasez\fi\noindent
    {\sezfont #1}\nobreak\vskip\spaziosottosez\nobreak\noindent}
%

\font\sezfont=cmbx10

\def\sezpreset#1{\global\sect@num=#1
    \immediate\write16{ !!! sez-preset = #1 }   }

\def\spaziosoprasez{50pt plus 60pt}
\def\spaziosottosez{15pt}
\def\spaziotitsez{5truemm}

\def\sref#1{\se@indefinito{@s@#1}\immediate\write16{ ??? \string\sref{#1}
    non definita !!!}
    \expandafter\xdef\csname @s@#1\endcsname{??}\fi\csname @s@#1\endcsname}

\def\autosez#1#2\par{
    \global\advance\sect@num by 1\if@ndoppia\global\eq@num=0\fi
    \xdef\la@sezionecorrente{\the\sect@num}
    \def\usa@getta{1}\se@indefinito{@s@#1}\def\usa@getta{2}\fi
    \expandafter\ifx\csname @s@#1\endcsname\la@sezionecorrente\def
    \usa@getta{2}\fi
    \ifodd\usa@getta\immediate\write16
      { ??? possibili riferimenti errati a \string\sref{#1} !!!}\fi
    \expandafter\xdef\csname @s@#1\endcsname{\la@sezionecorrente}
    \immediate\write16{\la@sezionecorrente. #2}
    \if@simboli
      \immediate\write\file@simboli{ }\immediate\write\file@simboli{ }
      \immediate\write\file@simboli{  Sezione 
                                  \la@sezionecorrente :   sref.   #1}
      \immediate\write\file@simboli{ } \fi
    \if@riferimenti
      \immediate\write\file@ausiliario{\string\expandafter\string\edef
      \string\csname\b@lank @s@#1\string\endcsname{\la@sezionecorrente}}\fi
    \goodbreak\vskip\spaziosoprasez
    \noindent\if@bozze\llap{\tt@bozze#1\ }\fi
      {\sezfont\the\sect@num.\hskip\spaziotitsez #2}\par\nobreak
    \vskip\spaziosottosez\nobreak\noindent}

\def\semiautosez#1#2\par{
    \gdef\la@sezionecorrente{#1}\if@ndoppia\global\eq@num=0\fi
    \if@simboli
      \immediate\write\file@simboli{ }\immediate\write\file@simboli{ }
      \immediate\write\file@simboli{  Sezione ** : sref.
          \expandafter\spo@glia\meaning\la@sezionecorrente}
      \immediate\write\file@simboli{ }\fi
    \s@ection#2\par}


\def\eqpreset#1{\global\eq@num=#1
     \immediate\write16{ !!! eq-preset = #1 }     }

\def\eqref#1{\se@indefinito{@eq@#1}
    \immediate\write16{ ??? \string\eqref{#1} non definita !!!}
    \expandafter\xdef\csname @eq@#1\endcsname{??}
    \fi\csname @eq@#1\endcsname}

\def\eqlabel#1{\global\advance\eq@num by 1
    \if@ndoppia\xdef\il@numero{\la@sezionecorrente.\the\eq@num}
       \else\xdef\il@numero{\the\eq@num}\fi
    \def\usa@getta{1}\se@indefinito{@eq@#1}\def\usa@getta{2}\fi
    \expandafter\ifx\csname @eq@#1\endcsname\il@numero\def\usa@getta{2}\fi
    \ifodd\usa@getta\immediate\write16
       { ??? possibili riferimenti errati a \string\eqref{#1} !!!}\fi
    \expandafter\xdef\csname @eq@#1\endcsname{\il@numero}
    \if@ndoppia
       \def\usa@getta{\expandafter\spo@glia\meaning
       \la@sezionecorrente.\the\eq@num}
       \else\def\usa@getta{\the\eq@num}\fi
    \if@simboli
       \immediate\write\file@simboli{  Equazione 
            \usa@getta :  eqref.   #1}\fi
    \if@riferimenti
       \immediate\write\file@ausiliario{\string\expandafter\string\edef
       \string\csname\b@lank @eq@#1\string\endcsname{\usa@getta}}\fi}

\def\autoreqno#1{\eqlabel{#1}\eqno(\csname @eq@#1\endcsname)
       \if@bozze\rlap{\tt@bozze\ #1}\fi}
\def\autoleqno#1{\eqlabel{#1}\leqno\if@bozze\llap{\tt@bozze#1\ }
       \fi(\csname @eq@#1\endcsname)}
\def\eqrefp#1{(\eqref{#1})}
\def\numeriadestra{\let\autoeqno=\autoreqno}
\def\numaeriasinistra{\let\autoeqno=\autoleqno}
\numeriadestra

\newcount\cit@num\global\cit@num=0

\newwrite\file@bibliografia
\newif\if@bibliografia
\@bibliografiafalse

\def\lp@cite{[}
\def\rp@cite{]}
\def\trap@cite#1{\lp@cite #1\rp@cite}
\def\lp@bibl{[}
\def\rp@bibl{]}
\def\trap@bibl#1{\lp@bibl #1\rp@bibl}

\def\refe@renza#1{\if@bibliografia\immediate        
    \write\file@bibliografia{
    \string\item{\trap@bibl{\cref{#1}}}\string
    \bibl@ref{#1}\string\bibl@skip}\fi}

\def\ref@ridefinita#1{\if@bibliografia\immediate\write\file@bibliografia{ 
    \string\item{?? \trap@bibl{\cref{#1}}} ??? tentativo di ridefinire la 
      citazione #1 !!! \string\bibl@skip}\fi}

\def\bibl@ref#1{\se@indefinito{@ref@#1}\immediate
    \write16{ ??? biblitem #1 indefinito !!!}\expandafter\xdef
    \csname @ref@#1\endcsname{ ??}\fi\csname @ref@#1\endcsname}

\def\c@label#1{\global\advance\cit@num by 1\xdef            
   \la@citazione{\the\cit@num}\expandafter
   \xdef\csname @c@#1\endcsname{\la@citazione}}

\def\bibl@skip{\vskip 0truept}


\def\stileincite#1#2{\global\def\lp@cite{#1}\global
    \def\rp@cite{#2}}
\def\stileinbibl#1#2{\global\def\lp@bibl{#1}\global
    \def\rp@bibl{#2}}

\def\citpreset#1{\global\cit@num=#1
    \immediate\write16{ !!! cit-preset = #1 }    }

\def\autobibliografia{\global\@bibliografiatrue\immediate
    \write16{ !!! Genera il file \jobname.BIB}\immediate
    \openout\file@bibliografia=\jobname.bib}

\def\cref#1{\se@indefinito                  
   {@c@#1}\c@label{#1}\refe@renza{#1}\fi\csname @c@#1\endcsname}

\def\cite#1{\trap@cite{\cref{#1}}}                        
\def\upcite#1{$^{\,\trap@cite{\cref{#1}}}$}               
\def\ccite#1#2{\trap@cite{\cref{#1},\cref{#2}}}           
\def\upccite#1#2{$^{\,\trap@cite{\cref{#1},\cref{#2}}}$}  
\def\cccite#1#2#3{\trap@cite{\cref{#1},\cref{#2},\cref{#3}}}          
\def\upcccite#1#2#3{$^{\,\trap@cite{\cref{#1},\cref{#2},\cref{#3}}}$} 
\def\ncite#1#2{\trap@cite{\cref{#1}--\cref{#2}}}          
\def\upncite#1#2{$^{\,\trap@cite{\cref{#1}-\cref{#2}}}$}  

\def\clabel#1{\se@indefinito{@c@#1}\c@label           
    {#1}\refe@renza{#1}\else\c@label{#1}\ref@ridefinita{#1}\fi}

\def\biblskip#1{\def\bibl@skip{\vskip #1}}           

\def\insertbibliografia{\if@bibliografia             
    \immediate\write\file@bibliografia{ }
    \immediate\closeout\file@bibliografia
    \catcode`@=11\input\jobname.bib\catcode`@=12\fi}


\def\commento#1{\relax} 
\def\biblitem#1#2\par{\expandafter\xdef\csname @ref@#1\endcsname{#2}}


\catcode`@=12

\input nuovisimboli.tex

\magnification 1200





\numerazionedoppia
\autobibliografia

\def\IM{\hbox{\rm{Im}}\,}
\def\ID{{\bf 1}}
\def\RE{\hbox{\rm{Re}}\,}

\def\ddt{{{\rm d}\over{\rm d} t}}

\def\ddtt{{{\rm d}\over{\rm d}\tau}}
\def\ddz{{{\rm d}\over{\rm d} z}}

\def\eps{\varepsilon}
\def\arg{{\rm arg}}

\centerline{\bf PICARD AND CHAZY SOLUTIONS TO THE PAINLEVE' VI EQUATION}
\vskip 0.2 cm
\centerline{\bf Marta Mazzocco}

\vskip 0.2 cm
\centerline{\it International School of Advanced Studies, SISSA-ISAS, 
Trieste.}
\vskip 0.5 cm

\noindent{\bf Abstract.}\quad
I study the solutions of a particular family of Painlev\'e VI equations 
with the parameters $\beta=\gamma=0$, $\delta={1\over2}$ and 
$2\alpha=(2\mu-1)^2$, for $2\mu\in\interi$. I show that the case 
of half-integer $\mu$ is integrable and that the solutions are of two 
types: the so-called Picard solutions and the so-called Chazy 
solutions. I give explicit formulae for them and completely determine their
asymptotic behaviour near the singular points $0,1,\infty$ and their nonlinear
monodromy. I study the structure of analytic continuation of the solutions to 
the PVI$\mu$ equation for any $\mu$ such that $2\mu\in\interi$.  As an 
application, I classify all the algebraic  
solutions. For $\mu$ half-integer, I show that they are in one to one 
correspondence with regular polygons or star-polygons in the plane.
For $\mu$ integer, I show that all algebraic solutions belong to a 
one-parameter family of rational solutions.

\vskip 0.2 cm\noindent SISSA preprint no. 89/98/FM, 13 August 1998. 
Revised 17 December 1998. 

\vskip 0.5 cm\semiautosez{1}{\bf 1. Introduction.}

In this paper I study the following particular case of Painlev\'e VI 
equation (see [Pain], [Gamb]):
$$
\eqalign{
y_{xx}=&{1\over2}\left({1\over y}+{1\over y-1}+{1\over y-x}\right) y_x^2 -
\left({1\over x}+{1\over x-1}+{1\over y-x}\right)y_x\cr
&+{1\over2}{y(y-1)(y-x)\over x^2(x-1)^2}\left[(2\mu-1)^2+
{x(x-1)\over(y-x)^2}\right],\cr}\eqno{PVI\mu}
$$
in the complex variable $x$, for resonant values of the parameter $\mu$,
i.e. $2\mu\in\interi$.
In the first part (see Sections 2, 3, 4) I show that, for any half-integer 
$\mu$, the PVI$\mu$ equation
is integrable and compute the solutions in terms of known special functions. 
In particular, I completely describe the asymptotic behaviour around the 
critical points $0,1,\infty$ for every branch of all solutions, and their
non linear monodromy. I show that for any half-integer $\mu$, PVI$\mu$ admits 
a countable set of algebraic solutions. In second part (see Sections 5, 6, 7, 
8) I describe the structure of analytic continuation of the solutions of the 
PVI$\mu$ equation for any $2\mu\in\interi$ and show that the algebraic 
solutions of PVI$\mu$ with half-integer $\mu$ are in one to one 
correspondence with regular polygons or star-polygons in the plane.
For $\mu$ integer, I show that there are no algebraic solutions except
a one parameter family of rational solutions.

The fact that the PVI$\mu$ equation with $\mu={1\over2}$ is integrable and
admits an infinite set of algebraic solutions, was already known to Picard, 
see [Pic]. The {\it Picard solutions,}\/ of PVI$\mu$ with $\mu={1\over2}$
are described in Section 2. They have the form
$$
y(x;\nu_1,\nu_2)=\wp (\nu_1 \omega_1 + \nu_2 \omega_2;\omega_1,\omega_2) 
+ {x+1\over3} 
$$
where $\omega_{1,2}(x)$ are two linearly independent solutions of the 
following Hypergeometric equation
$$
x(1-x)\omega''(x)+(1-2 x) \omega'(x)-{1\over4} \omega(x)=0,
$$ 
$\nu_1,\nu_2$ are complex numbers such that $0\leq\RE\nu_i<2$ and 
$\wp(u;\omega_1,\omega_2)$ is the Weierstrass elliptic function with the 
half-periods $\omega_1$, $\omega_2$.
I show that all the other PVI$\mu$ equations with half-integer 
$\mu\neq{1\over2}$, have ``more'' solutions. Let me briefly explain what 
I mean. Let the {\it solutions of Picard type}\/ be the solutions of 
PVI$\mu$ with $\mu+{1\over2}\in\interi\backslash\{1\}$ which are images 
via birational canonical transformations of Picard solutions. 
I show that, while the Picard solutions exhaust all the possible solutions 
of PVI$\mu$ with $\mu={1\over2}$, the solutions of Picard type do not cover 
all the possible solutions of PVI$\mu$ for all the other half integer values 
of $\mu$, i.e. for $\mu+{1\over2}\in\interi\backslash\{1\}$. Indeed, there 
exists a one-parameter family of transcendental solutions of PVI$\mu$ with 
$\mu+{1\over2}\in\interi\backslash\{1\}$, the so-called 
{\it Chazy solutions,}\/ which are not of Picard type. I describe the Chazy 
solutions of PVI$\mu$ with $\mu=-{1\over2}$ in Section 3. They are a one 
parameter family $y(x;\nu)$ of the form
$$
y={{1\over8}\left\{\left[\nu \omega_2+\omega_1+ 
2 x (\nu  \omega_2'+\omega_1')\right ]^2 -
4 x (\nu  \omega_2'+\omega_1')^2\right\}^2\over
(\nu \omega_2+\omega_1)(\nu \omega_2'+\omega_1') 
[2 (x-1) (\nu \omega_2'+ \omega_1')+\nu \omega_2+\omega_1]
[\nu \omega_2+\omega_1 +2 x (\nu \omega_2'+\omega_1')]},
$$
where $\omega_{1,2}(x)$ are chosen as above, $\omega'_{1,2}(x)$ are their 
derivatives with respect to $x$ and $\nu$ is a complex parameter.
The set of Chazy and Picard type solutions covers all the possible solutions 
of PVI$\mu$ with any half-integer $\mu\neq{1\over2}$. 
I compute explicitly the asymptotic behaviour of Picard and Chazy solutions 
for any choice of the parameters $(\nu_1,\nu_2)$ and $\nu$ respectively (see 
Lemma 2 and Lemma 5). I show that structure of the nonlinear monodromy is 
given by the action of $\Gamma(2)$ on $(\nu_1,\nu_2)$ and $\nu$, i.e. given 
a branch $y(x;\nu_1,\nu_2)$ (resp. $y(x;\nu)$) of a Picard (resp. Chazy) 
solution, all the other branches of the same solutions are of the form
$y(x;\tilde\nu_1,\tilde\nu_2)$ (resp. $y(x;\tilde\nu)$) with 
$$
\pmatrix{\tilde\nu_1\cr\tilde\nu_2\cr}=\pmatrix{a&b\cr c&d\cr}\pmatrix{
\nu_1\cr \nu_2\cr},\qquad \tilde\nu={a \nu+b\over c\nu+d}.
$$
Concerning the algebraic solutions, I show that, for any half-integer 
$\mu$, they form a countable set. They coincide with all the Picard type
solutions with rational $(\nu_1,\nu_2)$.

One of the main tools to prove the above results are the symmetry 
transformations between solutions of Painlev\'e equations with different 
values of the parameters (see [Ok]). In Section 4, I prove that all the 
solutions of PVI$\mu$ equations with any half-integer $\mu$, 
$\mu\neq{1\over2}$, are transformed via birational canonical transformations 
to solutions of the case $\mu=-{1\over2}$ and that the birational canonical 
transformations mapping the case $\mu=-{1\over2}$ to the case $\mu={1\over2}$ 
diverge when applied to the Chazy solutions. This is the reason why Chazy 
solutions are lost in the case of $\mu={1\over2}$.

Evenif the nonlinear monodromy of the solutions of PVI$\mu$ with half-integer 
$\mu$ is completely described in the first part of this paper (see Theorems 1
and 3), it is interesting to study the structure of analytic continuation of 
the solutions by a geometric approach which allows us to parametrize the 
algebraic solutions in terms of regular polyhedra or star polyhedra in the 
plane and to describe also the case of integer $\mu$. 
The main tool to study the structure of analytic continuation of the 
solutions to PVI$\mu$ with any $2\mu\in\interi$ is the isomonodromy 
deformation method (see [FlN], [ItN]), i.e PVI$\mu$ is treated as 
isomonodromy deformation equation of the auxiliary Fuchsian system
$$
{{\rm d} Y\over{\rm d} z}=
\left({A_0\over z}+{A_1\over z-1}+{A_x\over z-x}\right) Y,
\autoeqno{in2}
$$
where $A_0,A_1,A_x$ are $2\times2$ nilpotent matrices and 
$A_0+A_1+A_x=-A_\infty$, with
$$
\eqalign{
& A_\infty=\pmatrix{\mu&0\cr 0&-\mu\cr},\quad\hbox{for}\quad \mu\neq 0\cr
& A_\infty=\pmatrix{0&1\cr 0&0\cr},\quad\hbox{for}\quad \mu= 0\cr}.
$$
The technique is similar to the one developed in [DM], with some subtleties
due to the resonance of the matrix $A_\infty$, so I give a brief resume of it 
in Sections 5, 6. In Section 7, I give the explicit relation between 
algebraic solutions and monodromy data of the system 
\eqrefp{in2}. In Section 8, I show that the algebraic 
solutions of PVI$\mu$ with half-integer $\mu$ are in one to one 
correspondence with regular polygons or star-polygons in the plane and thus, 
the non linear monodromy of the algebraic elliptic curves of Weierstrass 
is described by the regular polygons or star-polygons in the plane.
For $\mu$ integer, I show that there are no algebraic solutions except
a one parameter family of rational solutions.
\vskip 0.3 cm

\noindent{\bf  Acknowledgments}\quad
I am indebted to B. Dubrovin who introduced me to the theory of Painlev\'e 
equations, constantly addressed my work and gave me lots of suggestions. 
I thank R. Conte for drawing the classical work of Picard (see [Pic])
to my attention and D. Guzzetti and E. Mukhin for their hints in Section 8. I 
thank also C. Reina and A. Zampa for many helpful discussions.

\semiautosez{2}{2. Picard Solutions.}
\vskip 0.3 cm

In this section, I describe the two parameter family of solutions of 
PVI$_{\mu={1\over2}}$ introduced by Picard (see [Pic]), their asymptotic 
behaviour and their monodromy. 

For the case $\mu=1/2$, Picard (see [Pic]) produced the following family of 
elliptic solutions:
$$
y(x)=\wp (\nu_1 \omega_1 + \nu_2 \omega_2;\omega_1,\omega_2) 
+ {x+1\over3} \autoeqno{P1}
$$
where $\omega_{1,2}(x)$ are two linearly independent solutions of the 
Hypergeometric equation
$$
x(1-x)\omega''(x)+(1-2 x) \omega'(x)-{1\over4} \omega(x)=0,
\autoeqno{C2}
$$
and $\nu_{1,2}$ are two complex numbers chosen such that 
$0\leq\RE \nu_{1,2}<2$ (all the other values of $\RE \nu_{1,2}$ can be 
reduced to this case by the periodicity of the Weierstrass function). 
I choose the branch cuts in the $x$-plane on the real axis $\pi_1=[-\infty,0]$
and $\pi_2=[1,+\infty]$.

\vskip 0.2 cm
\noindent{\bf Remark 1.}\quad Picard solutions obviously satisfy 
the Painlev\'e property. Indeed, they are regular functions of 
$\omega_{1,2}(x)$, which are analytic on the universal covering of
$\overline\complessi\backslash\{0,1,\infty\}$.

\proclaim Lemma 1. All the solutions of PVI$_{\mu=1/2}$ are of the form 
\eqrefp{P1}.

\noindent The proof is due to R. Fuchs (see [Fuchs] and [Man]).

\vskip 0.2 cm 
\noindent{\bf Remark 2.}\quad The general solution obtained by Hitchin 
(see [Hit]) in terms of Jacobi theta-functions in the case of 
PVI$_{{1\over 8},-{1\over 8},{1\over 8},{3\over 8}}$ can be obtained from
\eqrefp{P1} via the birational canonical transformation $w=w_1 w_2 w_1$ 
described in [Ok]. 

\vskip 0.3 cm
\noindent{\bf 2.1 Asymptotic behaviour and monodromy of the Picard 
solutions.}\quad
Here, I describe the monodromy of the Picard solutions and show that, 
for any $\nu_{1,2}\in \complessi^2\backslash\{(0,0)\}$, all 
the solutions \eqrefp{P1} have asymptotic behaviour of algebraic type (see 
(2.5) below).

First of all, let me fix a branch of a particular 
Picard solution \eqrefp{P1}, i.e. a pair of values $(\nu_1,\nu_2)$, and  
a branch of $\omega_{1,2}$ with respect to the above branch 
cuts $\pi_1$, $\pi_2$. For example, in a neighborhood of $0$, one can take:
$$
\omega^{(0)}_1(x)={\pi\over2} F\left({1\over2},{1\over2},1,x\right),
\qquad
\omega^{(0)}_2(x)={-i\over2} g\left({1\over2},{1\over2},1,x\right),
\autoeqno{C3.5}
$$
where
$$
F(a,b,c,z)=\sum_{k=0}^\infty {(a)_k(b)_k\over k!(c)_k} \, z^k,
$$
$$
g(a,b,z)=\sum_{k=0}^\infty {(a)_k(b)_k\over k!} \, z^k
[\ln z + \psi(a+k)+\psi(b+k)-2\psi(k+1)].
$$ 
Now, I fix some paths 
$\gamma_1$ and $\gamma_2$ along which the above basis is analytically 
continued to $1$ and to $\infty$ as in figure 1.

\midinsert
\centerline{\psfig{file=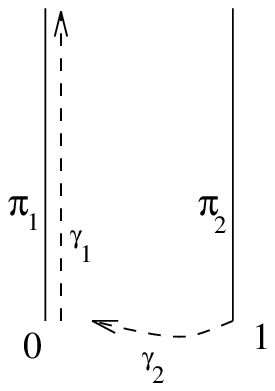,height=3cm}}\vskip 0.5 cm
\centerline{{\bf Fig.1.} The paths $\gamma_1$ and $\gamma_2$ along which the 
basis $\omega^{(0)}_{1,2}$ is analytically continued. }
\endinsert

Along the paths $\gamma_1$ and $\gamma_2$, the basis $\omega^{(0)}_{1,2}$ has 
the following analytic continuation:
$$
\eqalign{
\omega^{(0)}_1&\to\left\{\eqalign{
&\omega_1^{(1)}=-{1\over2}g\left({1\over2},{1\over2},1,1-x\right),
\qquad\hbox{as}\qquad x\to 1,\cr
&\omega_1^{(\infty)}={1\over2\sqrt{x}}
\left[i g\left({1\over2},{1\over2},1,{1\over x}\right)
+\pi F\left({1\over2},{1\over2},1,{1\over x}\right)\right],
\quad\hbox{as}\quad x\to \infty,\cr}\right.\cr
\omega^{(0)}_2&\to\left\{\eqalign{
&\omega_2^{(1)}={i\pi\over2}F\left({1\over2},{1\over2},1,1-x\right),
\qquad\hbox{as}\qquad x\to 1,\cr
&\omega_2^{(\infty)}=-{i\over2\sqrt{x}} 
g\left({1\over2},{1\over2},1,{1\over x}\right),
\qquad\hbox{as}\qquad x\to \infty.\cr}\right.\cr}
\autoeqno{P4}
$$

\proclaim Lemma 2. The above chosen branch of the solution \eqrefp{P1} of
PVI$\mu$ with $\mu={1\over2}$, for 
$\nu_{1,2}\in\complessi^2\backslash\{(0,0)\}$, $0\leq\RE \nu_{1,2}<2$, has 
the following asymptotic behaviour:
$$
y(x)\sim \left\{
\eqalign{
& a_0 x^{l_0} \left(1+{\cal O}(x^\eps)\right),
\qquad\hbox{as}\quad x\rightarrow 0,\cr
&1-a_1(1-x)^{l_1} \left(1+{\cal O}((1-x)^\eps)\right),
\qquad\hbox{as}\quad x\rightarrow 1,\cr
& a_\infty x^{1-l_\infty} \left(1+{\cal O}(x^{-\eps})\right),
\qquad\hbox{as}\quad x\rightarrow \infty,\cr}\right.
\autoeqno{P1.5}
$$
where $l_0,l_1,l_\infty$ are given by
$$
l_0=\left\{\eqalign{&\nu_2,\quad\hbox{if}\quad\RE\nu_2\leq 1\cr
&2-\nu_2,\quad\hbox{if}\quad\RE\nu_2>1\cr}\right.\quad
l_1=\left\{\eqalign{&\nu_1,\quad\hbox{if}\quad\RE\nu_1\leq1\cr
&2-\nu_1,\quad\hbox{if}\quad\RE\nu_1>1\cr}\right.
\autoeqno{P2.5}
$$
and
$$
l_\infty=\left\{\eqalign{&2+\nu_2-\nu_1,\quad\hbox{if}
\quad\RE(\nu_2-\nu_1)<0\cr
&\nu_2-\nu_1,\quad\hbox{if}
\quad\RE(\nu_2-\nu_1)\leq 1\cr
&2-(\nu_2-\nu_1),\quad\hbox{if}
\quad\RE(\nu_2-\nu_1)>1\cr}\right.
\autoeqno{P2.55}
$$
$a_0,a_1,a_\infty$ are three non-zero complex numbers depending on 
$\nu_{1,2}$ and $\eps>0$ is small enough. 

\vskip 0.2 cm
\noindent{\bf Remark 3.}\quad The solutions \eqrefp{P1}, for 
$(\nu_1,\nu_2)\in\complessi^2\backslash\{(0,0)\}$, with $\RE \nu_{1,2}>2$
or $\RE \nu_{1,2}<0$, can be reduced to the previous case thanks to the 
periodicity of the Weiestrass $\wp$ function. In this way, I show that the 
solutions \eqrefp{P1} have asymptotic behaviour of algebraic type for any 
$(\nu_1,\nu_2)\in \complessi^2\backslash\{(0,0)\}$, and that the exponents
$l_i$ have always real part in the interval $[0,1]$. 

\vskip 0.2 cm
\noindent Proof of Lemma 2. First, let me analyze the asymptotic 
behaviour of $y(x)$ as $x\to0$. Observe that, as $x\to0$ along any 
direction in the complex plane, the function $\tau(x)$, defined by 
$\tau={\omega_2\over\omega_1}$, has imaginary part that tends to infinity, 
while the real part remains limited. In fact, for $\omega_{1,2}(x)$ defined 
in \eqrefp{C3.5}:
$$
\tau=-i{g\left({1\over2},{1\over2},1,x\right)
\over \pi F\left({1\over2},{1\over2},1,x\right)}\sim
-i\pi\log|x| + {\arg(x)\over\pi},
$$
where, as $x\to0$, $\log|x|\to-\infty$, while $\arg(x)$ remains bounded 
for any fixed branch.
This fact permits to use the formula of the Fourier expansion of the 
Weierstrass function $\wp(u,\omega_1,\omega_2)$ (see [SG]):
$$
\wp(u,\omega_1,\omega_2)=-{\pi^2\over 12 \omega_1^2}+ 
{2\pi^2\over\omega_1^2} \sum_{k=1}^\infty {k q^{2 k}\over 1-q^{2k}}
\left(1-\cos{k \pi u\over\omega_1}\right) +
{\pi^2\over4\omega_1^2}\csc^2\left({\pi u\over2\omega_1}\right),
\autoeqno{P2}
$$
where $q=\exp({i\pi\omega_2\over\omega_1})$ and $u$ is such that
$$
-2\RE({\omega_2\over i\omega_1})<\RE({u\over i\omega_1})
<2\RE({\omega_2\over i\omega_1}).\autoeqno{P3}
$$
For $u=\nu_1\omega_1+\nu_2\omega_2$, \eqrefp{P3} reads:
$$
|\IM\nu_1 + \RE\nu_2 \IM\tau + \IM\nu_2 + \RE\tau|<2\IM\tau,
$$
that is always verified for $x\to0$ along any direction in the complex 
plane, and $|\RE\nu_2|<2$, because, $\IM\tau\to\infty$ while $\RE\tau$ 
remains limited. 

First, suppose that $\RE\nu_2\neq 0$. In this case,
$$
q=\exp\left({i \pi \omega_2\over\omega_1}\right)\qquad\hbox{with}\qquad
{i \pi \omega_2\over\omega_1}= {g\left({1\over2},{1\over2},1,x\right)\over
F\left({1\over2},{1\over2},1,x\right)}=\log(x) + 1 +{\cal O}(x).
$$
Then:
$$
\exp\left({i\pi u\over\omega_1}\right)= 
\exp(i\pi \nu_1)\exp\left(\nu_2\log(x)+\nu_2(1 +{\cal O}(x))\right)\sim
\exp(i\pi \nu_1+\nu_2) x^{\nu_2}(1 +{\cal O}(x)),
$$
and, for $\RE\nu_2>0$,
$$
\csc^2\left({\pi u\over2\omega_1}\right)
=-{4\over\exp\left({i \pi u\over\omega_1}\right)+
\exp\left(-{i \pi u\over\omega_1}\right)-2}
\sim-4 \exp(i\pi \nu_1) x^{\nu_2}+{\cal O}(x^{2 \nu_2}),
$$
$$
{\pi^2\over 12 \omega_1^2}\sim {1\over3}+{\cal O}(x),
$$
and
$$
{k q^{2 k}\over 1-q^{2k}} \left(1-\cos{k \pi u\over\omega_1}\right) 
\sim\left(k x^{2 k}+{\cal O}(x^{2k+1})\right) 
\left(-{\exp(-i\pi k \nu_1)\over2}x^{- k \nu_2} +{\cal O}(1)\right).
$$
As a consequence,
$$
y(x)\sim -4 \exp(i\pi \nu_1) x^{\nu_2} - 4 \exp(-i\pi \nu_1)
x^{2-\nu_2}+ {\cal O}\left(x^{3-\nu_2}\right)+ 
{\cal O}\left(x^{4-2\nu_2}\right)
+ {\cal O}\left(x^2\right)+ {\cal O}\left(x^{2 \nu_2}\right).
$$
This gives the required asymptotic behaviour around 0. 

For $\RE\nu_2=0$, with $\IM\nu_2\neq0$, 
$$
y(x)\sim {x\over 3}-
{4x^2(\alpha(x)-1)^4-\alpha^2\over\alpha(x)(\alpha(x)-1)}+ 
{\cal O}\left(x^2\right),
$$
where $\alpha(x)=\exp(i\pi \nu_1+\nu_2) x^{\nu_2}$ remains limited for 
$x\to0$ along any direction in the complex plane and
$$
-{4x^2(\alpha(x)-1)^4-\alpha^2\over\alpha(x)(\alpha(x)-1)}\sim
{\exp(i\pi \nu_1+\nu_2) x^{\nu_2}\over\exp(i\pi \nu_1+\nu_2) x^{\nu_2}-1},
$$
that gives the required asymptotic behaviour around 0 for $l_0=\nu_2$.
If $\nu_2=0$, then ${\pi u\over\omega_1}=\pi \nu_1$, and 
$$
y(x)\sim\csc^2\left({\pi \nu_1\over2}\right)+{\cal O}(x),
$$
which is again the required asymptotic behaviour around 0 for $l_0=0$.

The asymptotic behaviour of $y(x)$ as $x\to1$ can be obtained observing that,
as $x\to1$, the chosen branch of the basis $\omega_{1,2}$ is analytically
continued to $\omega^{(1)}_{1,2}$ given by \eqrefp{P4}. As a consequence
$$
y(x)\sim{x+1\over3}+
\wp\left(\nu_1\omega_1^{(1)}+\nu_2\omega_2^{(1)},
\omega_1^{(1)},\omega_2^{(1)}\right)
$$
and, with the change of variable $x=1-z$, one obtains, for $z\to0$,
$$
\eqalign{
y(z)\sim&{2-z\over3}+\wp\left(-\nu_1 \omega_2^{(0)}+\nu_2\omega_1^{(0)},
-\omega_2^{(0)},\omega_1^{(0)}\right)=\cr
=&{2-z\over3}+\wp\left((-\nu_2)\omega_1^{(0)}+\nu_1\omega_2^{(0)},
\omega_1^{(0)},\omega_2^{(0)}\right).\cr}
$$
Using the previously obtained asymptotic behaviour of $y(x)$ around zero, 
one immediately obtains the required asymptotic behaviour around $1$ with 
$l_1$ given by \eqrefp{P2.5}. 

Analogously, one can derive the asymptotic behaviour around $\infty$. The only
delicate point is that $\nu_2-\nu_1$ might be negative. In this case one takes
$2+\nu_2-\nu_1$ which is positive because $\nu_2-\nu_1>-2$.
This concludes the proof of the lemma. {\hfill QED}

\proclaim Theorem 1. The monodromy of the Picard solutions \eqrefp{P1} is 
described by the action of the group $\Gamma(2)$ on the parameters 
$(\nu_1,\nu_2)$, i.e. given a branch $y(x;\nu_1,\nu_2)$, all the other 
branches of the same solution are of the form $y(x;\tilde\nu_1,\tilde\nu_2)$ 
with all the $(\tilde\nu_1,\tilde\nu_2)$ such that
$$
\pmatrix{\tilde\nu_1\cr \tilde\nu_2\cr}=\pmatrix{a&b\cr c&d\cr}\pmatrix{
\nu_1\cr \nu_2\cr}\quad\hbox{for}\quad\pmatrix{a&b\cr c&d\cr}\in\Gamma(2).
$$

\vskip 0.2 cm 
\noindent Proof. \quad Let us fix a particular Picard solution \eqrefp{P1},
i.e. a particular pair of values of $(\nu_1,\nu_2)$. A branch is given by the 
choice of a branch of the basis $\omega_{1,2}$ of solutions of the 
hypergeometric equation \eqrefp{C2}.
As a consequence, the monodromy of the Picard solutions is described by the 
monodromy of the hypergeometric equation \eqrefp{C2}. This is given by the 
action of the group $\Gamma(2)$ on $\omega_{1,2}$. In fact, let us fix a 
basis $\gamma_0,\gamma_1$ of loops in 
$\pi_1\left(\overline\complessi\backslash\{0,1,\infty\}\right)$ like in 
figure 2. Let us consider $\omega_{1,2}$ chosen as in \eqrefp{C3.5}, 
\eqrefp{P4}. The result of the analytic continuation of $\omega_{1,2}^{(0)}$ 
along $\gamma_0$ is given by:
$$
\left(\eqalign{
&\omega_1^{(0)}\cr &\omega_2^{(0)}\cr}\right)\to
\pmatrix{
\omega_1^{(0)}\cr
2\omega_1^{(0)}+\omega_2^{(0)}\cr}=
\pmatrix{1&0\cr 2&1\cr}\left(\eqalign{
&\omega_1^{(0)}\cr &\omega_2^{(0)}\cr}\right)
$$
and the result of the analytic continuation of $\omega_{1,2}^{(1)}$ along 
$\gamma_1$ is given by:
$$
\left(\eqalign{
&\omega_1^{(1)}\cr &\omega_2^{(1)}\cr}\right)\to
\pmatrix{
\omega_1^{(1)}-2\omega_2^{(1)}\cr
\omega_2^{(1)}\cr}=
\pmatrix{1&-2\cr 0&1\cr}\left(\eqalign{
&\omega_1^{(1)}\cr &\omega_2^{(1)}\cr}\right).
$$
\midinsert
\centerline{\psfig{file=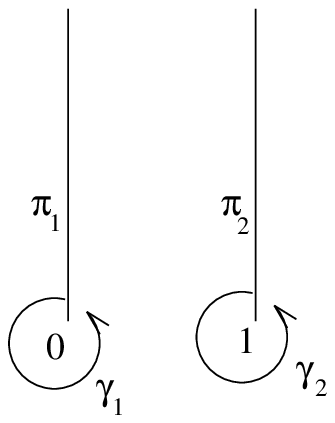,height=3cm}}\vskip 0.5 cm
\centerline{{\bf Fig.2.} The basis $\gamma_0$ and $\gamma_1$ of the loops
in $\pi_1\left(\overline\complessi\backslash\{0,1,\infty\}\right)$.}
\endinsert

The matrices $M_0=\pmatrix{1&0\cr 2&1\cr}$ and $M_1=\pmatrix{1&-2\cr 0&1\cr}$
generate the group $\Gamma(2)$. I stress that, as expected, 
$M_1 M_0=M_\infty^{-1}$, where $M_\infty$ gives the result of the 
analytic continuation of $\omega_{1,2}^{(\infty)}$ along 
$\gamma_0\cdot\gamma_1$, which is homotopic to $\gamma_\infty^{-1}$. 
Let $A=\pmatrix{a&b\cr c&d\cr}$ a matrix in $\Gamma(2)$. The new branch is 
given by $\tilde\omega_1=a\omega_1+b\omega_2$, 
$\tilde\omega_2=c\omega_1+d\omega_2$, and for any $\nu_1,\nu_2\in\complessi$ 
there exist $\tilde\nu_1,\tilde\nu_2$ such that 
$\nu_1\tilde\omega_1+\nu_2\tilde\omega_2=
\tilde\nu_1\omega_1+\tilde\nu_2\omega_2$. In fact
$\tilde\nu_1=a\nu_1+c\nu_2$ and $\tilde\nu_2=b\nu_1+d\nu_2$, i.e. 
$\tilde\nu_1,\tilde\nu_2$ are obtained from $\nu_1,\nu_2$
by the transformation induced by ${A^T}\in\Gamma(2)$.
This concludes the proof of theorem 1. {\hfill QED}

\vskip 0.3 cm
\noindent{\bf 2.2. Algebraic solutions to PVI$\mu$ with $\mu={1\over2}$.} 
\quad
Here I show that as known by Picard, PVI$\mu$ for $\mu={1\over2}$ admits a 
countable set of algebraic solutions and parameterize them. 

\vskip 0.2 cm
\noindent {\bf Definition.} A solution $y(x)$ is called {\it algebraic,}\/ if 
there exists a polynomial in two variables such that $F(x,y)\equiv0$.
An algebraic solution $y(x)$ is said {\it rational,}\/ if it is the ratio of 
two polynomials of variable $x$, with complex coefficients.

\proclaim Lemma 3. For any $\nu_1,\nu_2\in{\bf Q}$, the \eqrefp{P1} is a 
branch of an algebraic solution of PVI$_{\mu={1\over2}}$. Algebraic solutions 
form a countable set.  For any $\nu_1,\nu_2\in{\bf Q}$, the corresponding
algebraic solution contain the branch $y(x;{M\over N},0)$ where $(M,N)$ are 
coprime integers defined as follows: if $\nu_i={p_i\over q_i}$ for some pairs 
of coprime integers $(p_i,q_i)$, $i=1,2$, $N$ is the smallest common multiple 
of $q_1,q_2$ and $M$ is the largest common divisor of ${p_1 N\over q_1}$ and 
${p_2 N\over q_2}$.

\noindent Proof. For any $\nu_1,\nu_2\in{\bf Q}$, by the use of the 
addition and bisection formulae for the Weierstrass function, it is easy to 
see that $\wp(\nu_1\omega_1+\nu_2\omega_2,\omega_1,\omega_2)$ is an algebraic 
expression of the invariants $e_1,e_2,e_3$, which are given by 
$$
e_1=1-{x+1\over3},\qquad e_2=x-{x+1\over3},\qquad e_3 =-{x+1\over3}.
$$
This shows that for any $\nu_1,\nu_2\in{\bf Q}$, $y(x)$ is an algebraic 
function of $x$. It remains to show that two branches $y(x,\nu_1,\nu_2)$ and
$y(x,\tilde\nu_1,\tilde\nu_2)$ are branches of the same solution (up to the 
transformations $x\to1-x$, $y\to1-y$ and $x\to{1\over x}$, $y\to {y\over x}$) 
if and only if the correspondent integers $(M,N)$ and $(\tilde M,\tilde N)$, 
defined as in the statement of the theorem, coincide. Indeed, I show that
the ratio ${M\over N}$ is preserved under the analytic continuation. Due to
Theorem 1, the analytic continuation of a solution $y(x)$ is described by the 
action of $\Gamma(2)$ on $\nu_1,\nu_2$, that preserves the ratio ${M\over N}$.
Indeed, we can write $\nu_1=m_1 {M\over N}$ and $\nu_2=m_2{M\over N}$, where 
$m_{1,2}\in\interi$ and $(m_1,m_2)=1$, i.e. $m_1$ and $m_2$ are coprime 
integers. Consider a matrix $\pmatrix{a&b\cr c&d\cr}\in\Gamma(2)$. Then the 
new values of the parameters are given by 
$$
\pmatrix{\tilde\nu_1\cr\tilde\nu_2\cr}=
\pmatrix{a&b\cr c&d\cr}\pmatrix{\nu_1\cr\nu_2\cr}=
{M\over N}\pmatrix{a m_1+ b m_2\cr c m_1+ d m_2\cr}=
{M\over N}\pmatrix{\tilde m_1\cr \tilde m_2\cr},
$$
where $(\tilde m_1,\tilde m_2)=1$ because $a b- cd=1$, and then 
$(a,c)=(a,d)=(b,c)=(b,d)=1$. 

Now, consider any two numbers $\nu_1=m_1 {M\over N}$ and $\nu_2=m_2{M\over N}$,
for some given integers $(M,N,m_1,m_2)$ such that $(M,N)=1$ and $(m_1,m_2)=1$. 
There are three possibilities. i) $m_1$ and $m_2$ are odd integers. Then 
there exists a $\pmatrix{a&b\cr c&d\cr}\in\Gamma_2$ such that $a+b=m_1$, 
$c+d=m_2$. In fact, for any $\pmatrix{a&b\cr c&d\cr}\in\Gamma_2$, the numbers 
$a+b$ and $c+d$ are odd and coprime. As a consequence, the branch specified by 
$\nu_1=m_1{M\over N}$ and $\nu_2=m_2{M\over N}$, belongs to the same solution 
as the branch specified by $\nu_1={M\over N}=\nu_2$. ii) $m_1$ is even and 
$m_2$ is odd. Then there exists a $\pmatrix{a&b\cr c&d\cr}\in\Gamma_2$ such 
that $b=m_1$ and $d=m_2$.  As a consequence, the branch specified by 
$\nu_1= m_1{M\over N}$ and $\nu_2= m_2{M\over N}$, belongs to the same 
solution as the branch specified by $\nu_1=0$ and $\nu_2={M\over N}$.
ii) $m_1$ is odd and $m_2$ is even. Then there exists a 
$\pmatrix{a&b\cr c&d\cr}\in\Gamma_2$ such that $a=m_1$ and $c=m_2$.  
As a consequence, the branch specified by $\nu_1=m_1{M\over N}$ and 
$\nu_2=m_2{M\over N}$, belongs to the same solution as the branch 
specified by $\nu_1={M\over N}$ and $\nu_2=0$.
It easy to see that the above three cases are related one to the other by 
the transformations $x\to1-x$, $y\to1-y$ and $x\to{1\over x}$, $y\to 
{y\over x}$. This concludes the proof of the lemma. {\hfill QED}

\semiautosez{3}{3. Chazy Solutions.}

\vskip 0.7 cm
In this section, I introduce a one-parameter family of transcendental 
solutions of PVI$_\mu$, with $\mu=-{1\over2}$, compute their asymptotic 
behaviour and describe their monodromy.

\proclaim Theorem 2. There exists a one-parameter family of solutions of 
PVI$_\mu$ with $\mu=-{1\over2}$ of the form: 
$$
y(x)=
{{1\over8}\left\{\left[\nu \omega_2+\omega_1+ 
2 x (\nu  \omega_2'+\omega_1')\right ]^2 -
4 x (\nu  \omega_2'+\omega_1')^2\right\}^2\over
(\nu \omega_2+\omega_1)(\nu \omega_2'+\omega_1') 
[2 (x-1) (\nu \omega_2'+ \omega_1')+\nu \omega_2+\omega_1]
[\nu \omega_2+\omega_1 +2 x (\nu \omega_2'+\omega_1')]}
\autoeqno{C1}
$$
where $\omega_{1,2}(x)$ are two linearly independent solutions of the 
Hypergeometric equation \eqrefp{C2} and $\nu\in\complessi$ is the parameter.

\noindent Proof. Substitute \eqrefp{C1} in PVI$_{\mu=-{1\over2}}$. By
straightforward computations, it is easy to verify that if $\omega_{1,2}$ 
are solutions of \eqrefp{C2}, then \eqrefp{C1} 
is a solution of PVI$_{\mu=-{1\over2}}$ for any 
$\nu\in\complessi$. {\hfill QED} 

\vskip 0.2 cm 
Observe that the solutions \eqrefp{C1} are regular functions of 
$\omega_{1,2}(x)$, which are analytic on the universal covering of
$\overline\complessi\backslash\{0,1,\infty\}$, so they obviously satisfy 
the Painlev\'e property. I call {\it Chazy solutions}\/ the solutions 
\eqrefp{C1} and all the ones obtained from them via the transformations 
$x\to1-x$, $y\to1-y$ and $x\to{1\over x}$, $y\to{y\over x}$, which preserve 
the PVI$\mu$ equation. The reason of this name is that they correspond to the
following solution of WDVV equations in the variables $(t^1,t^2,t^3)$ 
(see [Dub]):
$$
F={(t^1)^2t^3\over2}+ {t^1(t^2)^2\over2}-{(t^2)^4\over 16} \gamma(t^3)
\autoeqno{C3}
$$
where the function $\gamma(t^3)$ is a solution of the equation of Chazy 
(see [Cha]):
$$
\gamma'''=6\gamma\gamma''-9{\gamma'}^2.
$$

\vskip 0.3 cm
\noindent{\bf 3.1. Derivation of the Chazy solutions.}\quad
I briefly outline how to derive \eqrefp{C1} from \eqrefp{C3}. Using the 
procedure explained in appendix E of [Dub], it is possible to show that the 
solution $y(x)$ of PVI$_{\mu=-{1\over2}}$ correspondent to the solution 
\eqrefp{C3} of WDVV equations is given by:
$$
\eqalign{
y(\tau)=&{\left (w_2(\tau) w_3(\tau)- w_1(\tau) w_2(\tau) - 
w_1(\tau) w_3(\tau)\right)^2\over 4 w_1(\tau) w_2(\tau) w_3(\tau) 
\left (w_1(\tau)-w_3(\tau)\right)},\cr
x(\tau)&={w_2(\tau)-w_1(\tau)\over w_3(\tau)-w_1(\tau)}\cr}
\autoeqno{C4}
$$
where $\tau=t^3$ and $(w_1,w_2,w_3)$ are solutions of the Halphen 
system (see [Hal]):
$$
\eqalign{
\ddtt w_1=&-w_1(w_2+w_3) + w_2 w_3,\cr
\ddtt w_2=&-w_2(w_1+w_3) + w_1 w_3,\cr
\ddtt w_3=&-w_3(w_1+w_2) + w_1 w_2,\cr} 
\autoeqno{C5}
$$
that is related to the Chazy equation. Indeed $(w_1,w_2,w_3)$ are the 
roots of the following cubic equation
$$
w^3+{3\over2}\gamma(\tau) w^2 +{3\over2} \gamma'(\tau) w +
{1\over4} \gamma''(\tau)=0.
$$
I want to derive \eqrefp{C1} from \eqrefp{C4}. The following lemma will 
be useful:

\proclaim Lemma 4. The transformation property
$$
\tilde w_i(\tau)={1\over c\tau +d} 
w_i\left({a\tau+b\over c\tau+d} \right) + {c\over c\tau+d},\qquad
\pmatrix{a&b\cr c&d\cr}\in PSL(2,\complessi),
$$
and the formulae
$$
w_1=-{1\over 2}\ddt\log{\lambda'\over\lambda},\quad
w_2=-{1\over 2}\ddt\log{\lambda'\over\lambda-1},\quad
w_3=-{1\over 2}\ddt\log{\lambda'\over\lambda(\lambda-1)},
\autoeqno{C6}
$$
where $\lambda(\tau)$ is a solution of the Schwartzian ODE:
$$
\{\tau,\lambda\} = {1\over2}\left[{1\over\lambda^2}+ {1\over(1-\lambda)^2}+
{1\over\lambda(1-\lambda)}\right],\autoeqno{C7}
$$ 
with
$$
\{\tau,\lambda\} = -\left[{\lambda'''\over\lambda'}-
{3\over2}\left({\lambda''\over\lambda'}\right)^2\right]
{1\over {\lambda'}^2},
$$ 
provide the general solution of \eqrefp{C5}.

\noindent The proof of this result can be found in [Tak].
\vskip 0.2 cm

The Schwartzian differential equation \eqrefp{C7} can be reduced to the 
hypergeometric equation \eqrefp{C2} via a standard procedure (see [Ince]). 
Let $\tau(\lambda)={\omega_1(\lambda)\over\omega_2(\lambda)}$. Then 
$\omega_{1,2}$ are two linearly independent solutions of the ODE:
$$
\omega''+p(\lambda) \omega'+q(\lambda)\omega=0, \autoeqno{C8}
$$
were $p(\lambda)$ and $q(\lambda)$ are two rational functions of $\lambda$
such that:
$$
2 q(\lambda)-{1\over2} p(\lambda)^2-p'(\lambda)=
{1\over2}\left[{1\over\lambda^2}+ {1\over(1-\lambda)^2}+
{1\over\lambda(1-\lambda)}\right].
$$
Requiring that \eqrefp{C8} is a hypergeometric equation, we obtain
$$
p(\lambda)={(1-2\lambda)\over(1-\lambda)\lambda},\qquad
q(\lambda)=-{1\over 4(1-\lambda)\lambda}.
$$
Using the formula \eqrefp{C4} for $x(\tau)$, and the formulae for the
Halphen functions \eqrefp{C6}, we see that $x(\tau)\equiv\lambda(\tau)$. 
As a consequence, \eqrefp{C7} is reduced to \eqrefp{C8} that coincides  
with \eqrefp{C2}. 

Putting $\tau(x)={\omega_1(x)\over\omega_2(x)}$, we can compute all the
derivatives $x'(\tau)$ and $x''(\tau)$ in terms of $x$, and by \eqrefp{C6},
$w_i(x)$. In this way we obtain \eqrefp{C1} by substitution in \eqrefp{C4}.

\vskip 0.3 cm
\noindent{\bf 3.2. Asymptotic behaviour and monodromy of the Chazy solutions.}
\quad Here I compute the asymptotic behaviour of the Chazy solutions and show 
that they are transcendental functions. Moreover I describe their nonlinear 
monodromy.

\proclaim Lemma 5. The solutions \eqrefp{C1}, for any $\nu\in\complessi$, 
and with branch cuts in the $x$-plane $\pi_1,\pi_2$ on the real axis, 
$\pi_1=[-\infty,0]$, $\pi_2=[1,-\infty]$, have the following asymptotic 
behaviour around the singular points $0,1,\infty$:
$$
y(x)\sim\left\{\eqalign{
&-\log(x)^{-2}  + b_0 \log(x)^{-3} +
{\cal O}\left(\log(x)^{-4} \right)
\quad\hbox{as}\quad x\to 0,\cr
&1+\log(1-x)^{-2} + b_1\log(1-x)^{-3} +
{\cal O}\left(\log(1-x)^{-4}\right),
\quad\hbox{as}\quad x\to 1\cr
&-x\log\left({1\over x}\right)^{-2} +
b_\infty x\log\left({1\over x}\right)^{-3}+ 
{\cal O}\left(\log\left({1\over x}\right)^{-4}\right),
\quad\hbox{as}\quad x\to \infty\cr}
\right.,\autoeqno{P4.5}
$$
where
$$
b_0= 1+{i\pi\over\nu}-4\log 2,\quad
b_1 = 2 [i\pi(\nu-1)-1+4\log 2],\quad
b_\infty=2[(\nu-1)(1-4 \log2)+i\pi].
$$

\noindent Proof. First of all one fixes a particular Chazy solution 
\eqrefp{C1}, i.e. a value $\nu$, and take a branch of it, i.e. a branch of 
$\omega_{1,2}$ for example \eqrefp{C3.5} and \eqrefp{P4}.
The correspondent branch $y(x)$ has the asymptotic behaviour \eqrefp{P4.5}
around the singular points $0,1,\infty$.{\hfill QED}.

\vskip 0.2 cm
Notice that the leading term of the asymptotic behaviour does not depend on 
the chosen particular solution, i.e. it does not depend on $\nu$. The 
dependence on $\nu$ appears in the second term. To derive the asymptotic 
behaviour of any other branch of $y(x)$, one can use the following: 

\proclaim Theorem 3. The monodromy of the Chazy solutions \eqrefp{C1} is 
described by the by the action of the group $\Gamma(2)$ on the
parameter $\nu$, for a fixed basis $\omega_{1,2}$, i.e. given 
a branch $y(x;\nu)$, all the other branches of the same solutions are of the 
form $y(x;\tilde\nu)$) with all $\tilde\nu$ such that
$$
\tilde\nu={a \nu+b\over c\nu+d}\quad\hbox{for}\quad
\pmatrix{a&b\cr c&d\cr}\in\Gamma(2).
$$

\vskip 0.2 cm 
\noindent The proof is analogous to the one of Theorem 1. \quad

\semiautosez{4}{4. Relations between half-integer values of $\mu$.}
\vskip 0.3 cm

We have seen in Section 2 that the PVI$_{\mu}$ for $\mu={1\over2}$ equation 
is integrable and its solutions are called {\it Picard solutions.}\/ 
I call {\it solutions of Picard type}\/ the solutions of PVI$\mu$ 
with any half-integer $\mu$, which are images of Picard solutions via the
birational canonical transformations described below. In this Section, I 
prove that, roughly speaking, while the Picard solutions exhaust all the 
possible solutions of PVI$_{\mu}$ for $\mu={1\over2}$, the Picard type ones 
do not exhaust all the possible solutions of the PVI$\mu$ equation with 
half-integer $\mu$, $\mu\neq{1\over2}$. All the PVI$\mu$ equations with any 
half-integer $\mu$, 
$\mu\neq{1\over2}$, are equivalent via birational canonical transfromations
to the case $\mu=-{1\over2}$, for which Picard type solutions and Chazy 
solutions are distinct and provide a complete set of solutions.
I recall that the so-called {\it singular solutions}\/ of the PVI$\mu$ equation
are $y(x)=0,1,\infty$. They are not solutions of the PVI$\mu$ equation because
the PVI$\mu$ equation itself becomes singular on them.

\proclaim Theorem 4. i) All the solutions $y(x)$ of PVI$\mu$ equations with 
$\mu+{1\over2}\in\interi$, $\mu\neq{1\over2}$, are mapped via birational  
canonical transformations to solutions of PVI$_{\mu=-1/2}$. ii) Chazy 
solutions of PVI$\mu$ are not of Picard type and vice-versa. iii) Chazy 
solutions and Picard type solutions exhaust all the possible solutions of 
PVI$\mu$ with $\mu+{1\over2}\in\interi$, $\mu\neq{1\over2}$.

\vskip 0.2 cm 
\noindent The proof is based on the following lemma:

\proclaim Lemma 6. The formula:
$$
\tilde y = y {\left(p_0 (y')^2+p_1 y'+p_2\right)^2
\over q_0(y')^4+ q_1(y')^3+
 q_2(y')^2+ q_3 y'+ q_4},\autoeqno{S1}
$$
where:
$$
\eqalign{
p_0&=x^2(x-1)^2,\cr
p_1& = 2 x(x-1)(y-1)[2\mu(y-x)-y],\cr
p_2& =y(y-1)[y(y-1)-4\mu(y-1)(y-x)+4\mu^2(y-x)(y-x-1)],\cr
q_0& =x^4(x-1)^4,\cr
q_1& =-4 x^3(x-1)^3 y(y-1),\cr
q_2& =2 x^2(x-1)^2 y(y-1)[3y(y-1)+4\mu^2(y-x)(1+x-3y)],\cr
q_3& =4x(x-1)y^2(y-1)^2[-y(y-1)-16\mu^3(y-x)^2+4\mu^2(y-x)(3y-x-1)],\cr
q_4& =y^2(y-1)^2\big\{y^2(y-1)^2+64\mu^3y(y-1)(y-x)^2,
-8\mu^2y(y-1)(y-x)(3y-x-1)+\cr
&+16\mu^4(y-x)^2[(x-1)^2+y(2+2x-3y)]\big\},\cr}
$$
transforms solutions of PVI$\mu$ to solutions of PVI$(-\mu)$, or 
equivalently of PVI$(\mu+1)$. 

\vskip 0.2 cm
\noindent Proof. Suppose $y(x)$ is a solution of PVI$\mu$. Take $\tilde y$ 
as in \eqrefp{S1} and compute $\tilde y'$ and $\tilde y''$ in terms of 
$(x,y,y')$ by derivating \eqrefp{S1} and substituting to $y''$ the righthand 
side of PVI$\mu$. Substituting $\tilde y$, $\tilde y'$ and $\tilde y''$ in 
PVI$_{-\mu}$ it is easy to verify that PVI$_{-\mu}$ is satisfied.

\vskip 0.2 cm
\noindent Proof of theorem 4. Let us consider a solution $y(x)$ of 
PVI$_\mu$ for any half-integer $\mu$.
We want to apply the transformation \eqrefp{S1} to $y(x)$. This is possible
if the denominator $Q(y_x,y,x,\mu)$ of the formula \eqrefp{S1} does not 
vanish identically on $y(x)$. In [DM], it is shown that $Q$ can identically 
vanish on solutions of PVI$\mu$ only for $\mu=-1/2$ or $\mu=0$. Let me 
consider the case PVI$_{\mu=-1/2}$. This is the same as PVI$_{\mu=3/2}$ 
(in fact $\mu$ and $1-\mu$ give the same 
value of the parameter $\alpha$ in PVI). For $\mu=3/2$ the denominator 
$Q(y_x,y,x,\mu)$ never vanishes and we can apply the transformation 
\eqrefp{S1}. Moreover $Q(y_x,y,x,\mu)$ never vanishes for any $\mu=3/2+n$,
$n\geq0$, so we can apply \eqrefp{S1} iteratively.
In this way we obtain all the PVI$_{\mu=\pm(3/2+n)}$ for any 
$n\geq0$. The above transformations are all invertible. In fact, starting 
from any PVI$_{\mu=\pm(3/2+n)}$, for $n\geq0$, we arrive at 
PVI$_{\mu=-1/2}$, via rational transformations of the form \eqrefp{S1}, the
determinant of which never vanishes. The idea of what happens it 
is shown in figure 2.

\midinsert
\centerline{\psfig{file=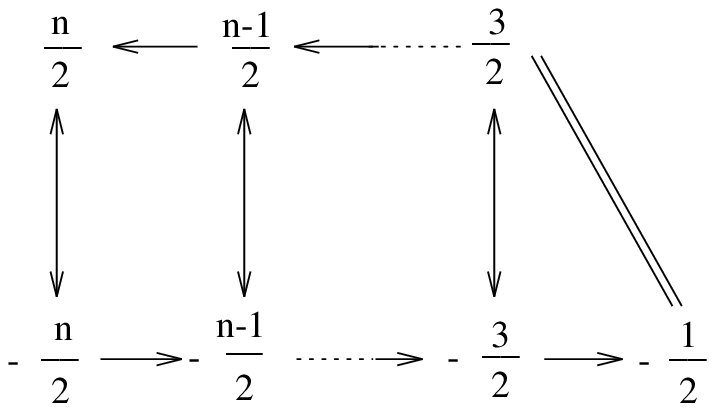,height=3cm}}\vskip 0.5 cm
\centerline{{\bf Fig.3.} The birational canonical transformations relating} 
\centerline{solutions of PVI$\mu$ equations with half-integer values of 
$\mu\neq{1\over2}$.}
\endinsert

We have proved the claim i) of Theorem 4. Claims ii) and iii) follow  
from the following two lemmas.

\proclaim Lemma 7. The one-parameter family of Chazy solutions exhaust all
the possible solutions of the differential equation $Q(y_x,y,x)\equiv0$.

\noindent Proof. Let us consider the algebraic 
differential equation $Q(y_x,y,x)\equiv0$. It has the following roots
$$
\eqalign{
y_x &=  {y(y-1)-\sqrt{y(y-1)}(y-x) + \sqrt{y(y-1)(y-x)}
\sqrt{2y-1 + 2 \sqrt{y(y-1)}}\over x(x-1)},\cr
y_x &= {y(y-1)-\sqrt{y(y-1)}(y-x) - \sqrt{y(y-1)(y-x)}
\sqrt{2y-1 + 2 \sqrt{y(y-1)}}\over x(x-1)},\cr
y_x &= {y(y-1)+\sqrt{y(y-1)}(y-x) + \sqrt{y(y-1)(y-x)}
\sqrt{2y-1 - 2 \sqrt{y(y-1)}}\over x(x-1)},\cr
y_x &= {y(y-1)+\sqrt{y(y-1)}(y-x) - \sqrt{y(y-1)(y-x)}
\sqrt{2y-1 - 2 \sqrt{y(y-1)}}\over x(x-1)}.\cr}
$$
All the above differential equations are equivalent: the first is mapped in
the third by the transformation $x\to 1-x$, $y\to 1-y$, to the forth by 
$x\to{1\over x}$, $y\to{y\over x}$ and to the second by $x\to{1\over 1-x}$,
$y\to{1-y\over 1-x}$. All these transformations preserve the class of 
Chazy solutions of PVI$_{\mu=-{1\over2}}$. Thus it is enough to prove that the 
one-parameter family of Chazy solutions exhaust all
the possible solutions of the first differential equation.
Indeed it easy to verify that $y(x,\nu)$ of the form \eqrefp{C1} solves 
it for any value of the parameter $\nu$. To conclude, we have to 
show that $\forall(x_0,y_0)\in\complessi\times\complessi$, there exists a 
value of the parameter $\nu$ such that
$$
y(x_0,\nu)=y_0.\autoeqno{C10}
$$        
If $y_0=\infty$, we can take $\nu=\nu_\infty$ such that 
$w'_1(x_0)+\nu_\infty w'_2(x_0)=0$. Let us suppose that $y_0\neq\infty$. 
Then $w'_1(x_0)+\nu w'_2(x_0)\neq0$ for every $\nu\neq \nu_\infty$ and 
$y(x)$ can be written in the form:
$$
y(x)={\left[(W(x,\nu)+ 2 x)^2-4 x\right]^2\over 8 W(x,\nu)
(W(x,\nu)+ 2x)[W(x,\nu)+ 2(x-1)]}
$$
where $W(x,\nu)={w_1(x)+\nu w_2(x)\over w'_1(x)+\nu w'_2(x)}$. If we show 
that given any $(x_0,W_0)$, there exists $\nu_0$ such that 
$W(x_0,\nu_0)=W_0$, we are done.  Indeed, for
$$
\nu_0=-{w_1(x_0) - W_0 w_1'(x_0) \over w_2(x_0)-W_0 w_2'(x_0)}
$$
$W(x_0,\nu_0)=W_0$.  {\hfill QED}

\proclaim Lemma 8. The denominator $Q(y_x,y,x)$ never vanishes on Picard 
type solutions.

\noindent Proof. Consider any Picard solution $y(x)$ and its
correspondent Picard type solution $\tilde y$, obtained by the 
transformation \eqrefp{S1}. I want to show that $\tilde y(x)$ is such that 
$Q(\tilde y_x,\tilde y,x)\neq0$. By straightforward computations, one obtains:
$$
\eqalign{
Q(\tilde y_x,\tilde y,x) &=
(x-1)^2x^2\left (y - y^2 - x y_x^2 + x^2 y_x^2\right)^4
\left (y^2-y - 2 x y_x(y - 1) - x y_x^2 + x^2 y_x^2\right)^4
\times\cr
&\times\left (y^2-y - 2 y y_x (x-1)  - x y_x^2 + x^2 y_x^2\right)^4\,
\bigg\{ y^2 (y-1)^2 - 4 y^2 (y-1)^2 y_x+\cr
& + 2 (y-1) yy_x^2 (4 x y-x -x^2 - 2 y)
-  4 (x-1) x ( y-1) y y_x^3+ (x-1)^2 x^2 y_x^4\bigg\}^{-1}.\cr}
$$
The above quantity can not vanish on any Picard solution $y$. In fact none 
of the polynomials
$$
\eqalign{
Q_1(y_x,y,x)=&\, y - y^2 - x y_x^2 + x^2 y_x^2,\cr
Q_2(y_x,y,x)=&\, y^2-y - 2 x y_x(y - 1) - x y_x^2 + x^2 y_x^2,\cr
Q_3(y_x,y,x)=&\, y^2-y - 2 y y_x (x-1)  - x y_x^2 + x^2 y_x^2,\cr}
$$
can vanish on any Picard solution. Indeed, eliminating $y_{xx}$ and $y_x$, 
form the system
$$
\eqalign{
&y_{xx}={1\over2}\left({1\over y}+{1\over y-1}+{1\over y-x}\right)
y_x^2 -\left({1\over x}+{1\over x-1}+{1\over y-x}\right)y_x
+{y(y-1)\over 2x(x-1)(y-x)},\cr
&Q_i(y_x,y,x,\mu)=0,\cr
&{{\rm d}\over {\rm d} x} Q_i(y_x,y,x,\mu)=0,\cr}
$$
for each $i=1,2,3$, we obtain the following resultants:
$$
( x-1) x (x - y)^2 ( y-1)^3 y^3,\quad    
(x-1)^3 x (x - y)^4 (y-1)^4 y^2,\quad
(x-1) x^3 (x - y)^4 ( y-1)^2 y^4,
$$
which never vanish for nosingular solutions of PVI$\mu$. This concludes the 
proof of Lemma 8. 
{\hfill QED}

\vskip 0.2 cm
Now, I conclude the proof of Theorem 4.
Claim ii) follows from the fact that, thanks to Lemma 5, $Q$ does not vanish 
on Picard type solutions, while, thanks to Lemma 7, it vanishes on Chazy 
solutions. Claim iii) is due to the fact that any solution of 
PVI$_{\mu={1\over2}}$ such that $Q(y_x,y,x,\mu)\neq0$ is necessarily of Picard
type. In fact, being $Q(y_x,y,x,\mu)\neq0$, the birational canonical 
transformation \eqrefp{S1} can be applied to $y$ and it gives rise to a 
Picard solution. {\hfill QED}

\vskip 0.2 cm
\noindent{\bf Remark 4.}\quad In the notations of the paper [Ok], the symmetry
from $\mu=-1/2$ to $\mu=1/2$ is given by the transformation $l_3^2$ applied
on the canonical coordinates $(p,q)$ of PVI$\mu=-{1\over2}$. The condition 
$Q(y_x,y,x)\equiv0$ is exactly the condition that the intermediate 
coordinates $l_3(p,q)$ are singular, i.e. the associated auxiliary 
Hamiltonian $h(l_3(p,q))$ is linear in $t$. In particular, this means
that the transformation $l_3$ can map non-singular solutions into singular 
ones. This is exactly what happens when we apply \eqrefp{S1} on the Chazy
solutions. 

\vskip 0.3 cm
\noindent{\bf 3.1. Chazy solutions as limit of Picard type solutions.}
\quad For $\nu_1=\nu_2=0$ the Weierstrass $\wp$-function has 
a pole and the correspondent function $y(x)$ defined by \eqrefp{P1} does
not exist.

\proclaim Lemma 9. Chazy solutions of PVI$_{\mu=-{1\over2}}$ can be obtained 
as the limit for $\nu_{1,2}\to0$, with ${\nu_2\over\nu_1}=\nu$, of the 
Picard type solution obtained applying the symmetry \eqrefp{S1} to the 
solution \eqrefp{P1} of PVI$\mu$ with $\mu={1\over2}$. 

The above result is not surprising, in fact, as observed above, Chazy 
solutions are transformed, via the symmetry \eqrefp{S1}, to singular 
solutions of PVI$_{\mu={1\over2}}$ which are identically equal to $\infty$.

\vskip 0.2 cm
\noindent Proof of Lemma 9. Consider a solution $y(x)$ of 
PVI$_{\mu={1\over2}}$, $y(x)$ given by the formula \eqrefp{P1}. Fix the 
ratio ${\nu_2\over\nu_1}=\nu$ and let $\nu_{1,2}\to \infty$.
Since the Weierstrass function has a pole of order two in $0$, one has
$$
\lim_{\nu_1\to0} \nu_1^2 y(x)={1\over(\omega_1+\nu\omega_2)^2},
$$
and
$$
\lim_{\nu_1\to0} \nu_1^2 y'(x)={1\over(\omega'_1+\nu\omega'_2)^2},
$$ 
and, applying the transformation \eqrefp{S1} for $\mu={1\over2}$ to 
$y(x)$ given by \eqrefp{P1}, and taking the limit as $\nu_{1,2}\to0$ with 
fixed ratio ${\nu_2\over\nu_1}=\nu$, one obtains the formula 
\eqrefp{C1}. This concludes the proof of the lemma.{\hfill QED}

\vskip 0.3 cm
\noindent{\bf 4.2. Algebraic solutions.}\quad
Here, I classify all the algebraic solutions of PVI$\mu$, for any
$\mu+{1\over2}\in\interi$. As shown in Section 2.2, PVI$\mu$ with 
$\mu={1\over2}$ admits a countable set of algebraic solutions. Now I show 
that for all the other half integer values of $\mu$, the algebraic solutions 
are all of Picard type, so they are equivalent via birational canonical 
transformations to the ones of Lemma 3.

\proclaim Theorem 5. For any half-integer $\mu$, PVI$_\mu$ admits a 
countable family of algebraic solutions. All the algebraic solutions are 
of Picard type for some $\nu_1,\nu_2\in{\bf Q}$, $0\leq\nu_{1,2}<2$.

\noindent Proof. The algebraic solutions are preserved under the 
transformations \eqrefp{S1}. So we obtain a countable family of algebraic 
solutions PVI$_\mu$ for any half-integer $\mu$. Moreover Chazy solutions 
are transcendental, so the algebraic solutions can only be of Picard 
type. 

{\hfill QED}

\vskip 0.2 cm
For example, we can recover the solutions found in [Dub], (E.34a), (E.36) 
and (E.37). In fact they are mapped by the symmetry \eqrefp{S1}, 
respectively to
$$
y = {(s-1)^2\over(s-3)(1 + s)},\qquad
x  ={(s-1)^3(3 + s)\over(s-3)(1 + s)^3},\autoeqno{a2}
$$
that is a Picard solution with $N=3$ and $M=2$, to
$$
y ={2+s\over4},\qquad
x  ={(s+2)^2\over8 s}  \autoeqno{b2}
$$
that is a Picard solution with $N=2$, i.e. $y(x)= x+\sqrt{(x-1)x}$ and to
$$
y = {3 (3 - t) (1 +t)\over(3 + t)^2},\qquad 
x = {(3-t)^3(1+t)\over(1-t)(3+t)^3},\autoeqno{g2}
$$
that is a Picard solution with $N=3$, $M=1$.
In section 8, I show how these solutions correspond to the 
affine Weyl groups $A_2$, $B_2$ and $G_2$ respectively and show that all 
the other algebraic solutions correspond to different presentations of
the dihedral group. 




\semiautosez{5}{5. Painlev\'e VI equation as isomonodromy deformation 
equation.}\vskip 0.5 cm

First of all, I briefly describe how PVI$\mu$, for $2\mu\in\interi$, can  
be interpreted as the {\it isomonodromy deformation equation}\/ of the 
following auxiliary Fuchsian system (see [JMU], [ItN], [FlN]) with  
four regular singularities at $z=u_1,u_2,u_3,\infty$:
$$
\ddz Y={\cal A}(z)Y,
\qquad z\in\overline{\bf C}\backslash\{u_1,u_2,u_3,\infty\}
\autoeqno{M1}
$$
where
$$
{\cal A}(z)=
{{\cal A}_1\over z-u_1} + {{\cal A}_2\over z-u_2} + {{\cal A}_3\over z-u_3},
$$
$u_1,u_2,u_3$ being pairwise distinct complex numbers and ${\cal A}_i$ 
being $2\times 2$ matrices satisfying the following conditions:
$$
{\cal A}_i^2=0\quad\hbox{and}\quad
 -{\cal A}_1-{\cal A}_2-{\cal A}_3={\cal A}_\infty
\autoeqno{M2}
$$
where ${\cal A}_\infty=\pmatrix{\mu & 0\cr 0 & -\mu\cr}$ for $\mu\neq0$ and
${\cal A}_\infty=\pmatrix{0 & 1\cr 0 & 0\cr}$ for $\mu=0$. 

The solution $Y(z)$ of the system \eqrefp{M1} is a multi-valued analytic 
function on the punctured Riemann sphere, ${\bf C}\backslash\{u_1,u_2,u_3\}$, 
and its multivaluedness is 
described by the so-called {\it monodromy matrices.}\/  To define them, we 
fix the basis $\gamma_1,\gamma_2,\gamma_3$ of loops in the fundamental group
$$
\pi_1\left(\overline{\bf C}\backslash\{u_1,u_2,u_3,\infty\},\infty\right),
$$ 
as in figure 4 and a fundamental matrix.
\midinsert
\centerline{\psfig{file=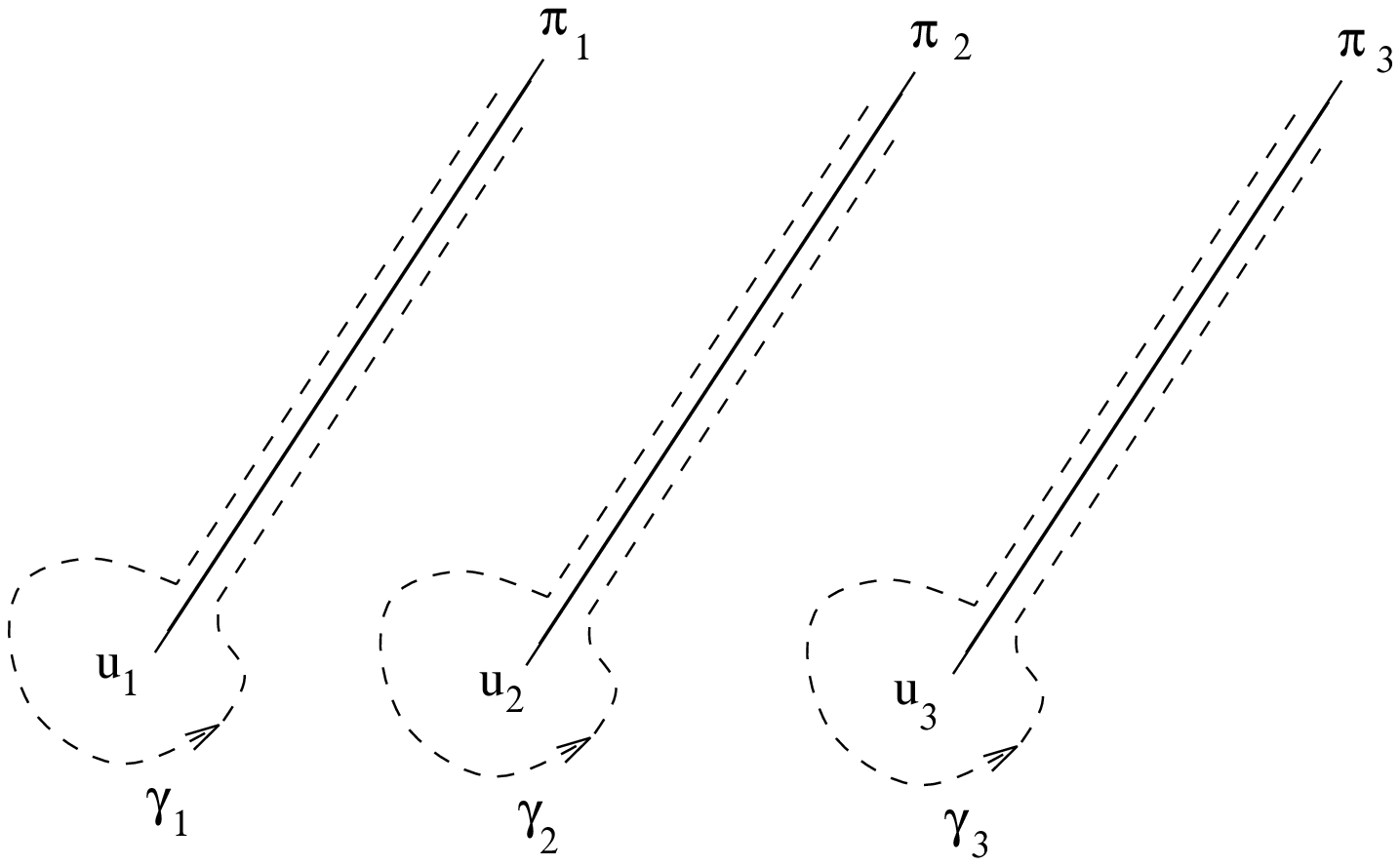,height=4cm}}\vskip 0.5 cm
\centerline{{\bf Fig.4}: The cuts $\pi_i$ between the singularities $u_i$, 
ordered according to the order of} 
\centerline{the points $u_1,u_2,u_3$, and the oriented loops 
$\gamma_i$, starting and finishing at infinity,}
\centerline{going around $u_i$ in positive direction ($\gamma_i$ is oriented 
counter-clockwise, $u_i$ lies}
\centerline{inside, while the other singular points lie outside) and not 
crossing the cuts $\pi_i$.}\vskip 0.1 cm
\endinsert

\proclaim Proposition 1. There exists a fundamental matrix of the system 
\eqrefp{M1} of the form
$$
Y_\infty=\left(\ID +{\cal O}({1\over z})\right) z^{-{\cal A}_\infty}
z^{R},\quad\hbox{as}\quad 
z\rightarrow\infty,
\autoeqno{M3}
$$
where the matrix $R$ is defined as follows for $2\mu=n\in\interi$,
$$
\eqalign{
&\hbox{for}\quad n>0,\quad 
R_{12}=\sum_{k=1}^3\left({\cal A}_k\right)_{12} u_k^n+
\sum_{l=1}^{n-1}\left(G^{(n-l)}
\sum_{k=1}^3 {\cal A}_k u_k^l \right)_{12},
\quad R_{11}=R_{21}=R_{22}=0,\cr
&\hbox{for}\quad n<0,\quad
R_{21}=\sum_{k=1}^3\left({\cal A}_k\right)_{21} u_k^n+
\sum_{l=1}^{n-1}\left(G^{(n-l)}
\sum_{k=1}^3 {\cal A}_k u_k^l\right)_{21},
\quad R_{11}=R_{12}=R_{22}=0,\cr
&\hbox{for}\quad n=0,\quad
R_{11}=R_{12}=R_{21}=R_{22}=0,\cr}
$$
where for $l=1,2,\cdots|n|-1$, $G^{(l)}$ are uniquelly determined by
$$
G^{(l)}_{11}=-\sum_{k=1}^3 \left({\cal A}_k\right)_{11} {u_k^l\over l}-
\sum_{i=1}^{l-1}\left(G^{(l-i)}\sum_{k=1}^3 {\cal A}_k u_k^i\right)_{11},
\autoeqno{DG1}
$$
$$
G^{(l)}_{12}=-\sum_{k=1}^3 \left({\cal A}_k\right)_{12} {u_k^l \over l-2\mu}-
\sum_{i=1}^{l-1}\left(G^{(l-i)}\sum_{k=1}^3 {\cal A}_k u_k^i\right)_{12},
\autoeqno{DG2}
$$
$$
G^{(l)}_{21}=-\sum_{k=1}^3 \left({\cal A}_k\right)_{21} {u_k^l\over l+2\mu}-
\sum_{i=1}^{l-1}\left(G^{(l-i)}\sum_{k=1}^3 {\cal A}_k u_k^i\right)_{21},
\autoeqno{DG3}
$$
$$
G^{(l)}_{22}=-\sum_{k=1}^3 \left({\cal A}_k\right)_{22} {u_k^l\over l}-
\sum_{i=1}^{l-1}\left(G^{(l-i)}\sum_{k=1}^3 {\cal A}_k u_k^i\right)_{22},
\autoeqno{DG4}
$$
and $z^\mu$ is defined as $e^{\mu\log z}$, with the choice of the principal 
branch of the logarithm with the branch-cut along the common direction of 
the cuts $\pi_1,\pi_2,\pi_3$ of figure 4.

\vskip 0.2 cm 
\noindent Proof. Consider the Fuchsian system near infinity. Perform a change 
of variable $z\to{1\over z}$, so that the Fuchsian system becomes
$$
\ddz Y=-{1\over z}\left({\cal A}_\infty + 
\sum_{l=1}^\infty {\cal A}^{(l)} z^l \right)Y,
$$
where
$$
{\cal A}^{(l)}= {\cal A}_1 u_1^l+ {\cal A}_2 u_2^l + {\cal A}_3 u_3^l.
$$
We look for a gauge 
$$
\tilde Y=G(z) Y,\qquad G(z)=\ID + \sum_{l=1}^\infty G^{(l)} z^l
$$
such that
$$
\ddz\tilde Y=-{1\over z}\left({\cal A}_\infty + 
\sum_{l=1}^N R^{(l)} z^l \right)Y,
$$
where 
$$
R^{(l)}_{ij}={\cal A}^{(l)}_{ij}+
G^{(l)}_{ij}(l+({\cal A}_\infty)_{jj}-({\cal A}_\infty)_{ii})
+\sum_{i=1}^{l-1}\left(G^{(l-i)}{\cal A}^{(i)}\right)_{ij}.
$$
One chooses $R^{(l)}_{ij}\neq 0$ if and only if 
$({\cal A}_\infty)_{ii}-({\cal A}_\infty)_{jj}=l$. Since
$({\cal A}_\infty)_{11}=\mu$ and $({\cal A}_\infty)_{22}=-\mu$, there 
exists a unique $l$ such that the above condition is fullfilled, i.e. 
$l=N=|2\mu|$. All the other $R^{(l)}$ are zero and for every $l<N$, 
$G^{(l)}$ is uniquelly determined by
\eqrefp{DG1}, \eqrefp{DG2}, \eqrefp{DG3}, \eqrefp{DG4}. {\hfill QED}

\vskip 0.2 cm
The monodromy matrix $M_\gamma$, by definition, is a constant invertible 
$2\times 2$ matrix such that 
$$
\gamma[Y_\infty(z)]= Y_\infty(z) M_\gamma,
$$
where $\gamma[Y_\infty(z)]$ is the analytic continuation of $Y_\infty(z)$ 
along the loop $\gamma$. Particularly, the matrix 
$M_\infty:=M_{\gamma_\infty}$ is given by:
$$
M_\infty=\exp{2\pi i({\cal A}_\infty+R)},\autoeqno{N9}
$$
where $\gamma_\infty$ is a loop around infinity in the clock-wise direction.
The matrices $M_i:=M_{\gamma_i}$, where the $\gamma_i$ for $i=1,2,3$ 
are the generators of the fundamental group, are given by:
$$
M_i= C_i^{-1} \exp(2\pi i J) C_i,\qquad i=1,2,3,\autoeqno{N4}
$$
where the matrices $C_1,C_2, C_3$ are  the so-called {\it connection 
matrices.}\/ The matrices $M_1$, $M_2$, $M_3$ generate the {\it monodromy 
group of the system,}\/ and satisfy the following relations:
$$
\det(M_i)=1, \quad {\rm Tr}(M_i)=2,\quad\hbox{for}\quad i=1,2,3,
\autoeqno{N5}
$$
with $M_i=\ID$ if and only if ${\cal A}_i=0$. Since the loop 
$(\gamma_1 \gamma_2 \gamma_3)^{-1}$ is homo-topic to $\gamma_\infty$, the 
following relation holds true:
$$
M_\infty M_3 M_2 M_1=\ID.\autoeqno{N6}
$$

\proclaim Lemma 10. For $\mu\neq 0$, under the assumption that $R\neq0$, two 
Fuchsian systems of the form \eqrefp{M1} with the same poles $u_1$, $u_2$ and 
$u_3$, and the same resonant value of $\mu$, coincide if and only if they 
have the same monodromy matrices $M_1$, $M_2$, $M_3$, with respect to the 
same basis of the loops $\gamma_1$, $\gamma_2$ and $\gamma_3$ and the same 
value of $R$. For $\mu=0$ two Fuchsian systems of the form \eqrefp{M1} with 
the same poles $u_1$, $u_2$ and $u_3$ coincide if and 
only if they have the same monodromy matrices $M_1$, $M_2$, $M_3$, with 
respect to the same basis of the loops $\gamma_1$, $\gamma_2$ and $\gamma_3$.

\noindent Proof. Fix for example $\mu>0$, i.e. $M_\infty$ 
upper-triangular. Suppose that there are two Fuchsian systems of the form 
\eqrefp{M1} with the same poles $u_1$, $u_2$ and $u_3$, and the same value 
of $\mu$, the same monodromy matrices $M_1$, $M_2$, $M_3$, the same value 
of $R$. The fundamental matrices at $\infty$ of the form \eqrefp{M3} exist. 
All the fundamental matrices of the form
$$
Y_\infty=\left(\ID +{\cal O}({1\over z})\right) z^{-{\cal A}_\infty}
z^{R}\pmatrix{1&a\cr0&1\cr}
$$
with any $a\in\complessi$, give the same monodromy matrix at infinity. So, 
for a given $M_\infty$, the fundamental matrices at $\infty$ of the two 
Fuchsian system can be fixed as
$$
Y_\infty^{(i)}=\left(\ID +{\cal O}({1\over z})\right) z^{-{\cal A}_\infty}
z^{R} \pmatrix{1&a^{(i)}\cr0&1\cr},\quad i=1,2
$$
for some $a^{(1)}$ and $a^{(2)}$. For any choice of $a^{(1)}$ and $a^{(2)}$, 
the following matrix:
$$
Y(z):= Y_\infty^{(2)}(z)Y_\infty^{(1)}(z)^{-1}.
$$
$Y(z)$ is an analytic function around infinity: 
$$
Y(z)=1+{\cal O}\left({1\over z}\right),\quad\hbox{as}\, z\rightarrow\infty,
$$
and, since the monodromy matrices coincide, $Y(z)$ is a single valued function 
on the punctured Riemann sphere $\overline\complessi\backslash\{u_1,u_2,u_3\}$.
As shown in [DM] Lemma 4.1, near the point $u_i$, we can choose the 
fundamental matrices $Y_i^{(1)}(z)$ and $Y_i^{(2)}(z)$ in such a way that
$$
Y_\infty^{(1),(2)}(z)= Y_i^{(1),(2)}(z)C_i \quad i=1,2,3.
$$
with the same connection matrices $C_i$ and
$$
Y_i=G_i\left(\ID +{\cal O}(z-u_i)\right)(z-u_i)^{J}.
$$
Then near the point $u_i$, $Y(z)$ is again analytic
$$
Y(z)=G_i^{(2)}\left(\ID+{\cal O}(z-u_i) \right)
\left[G_i^{(1)}\left(\ID+{\cal O}(z-u_i) \right)\right]^{-1}.
$$
This proves that $Y(z)$ is an analytic function on all $\overline\complessi$ 
and then, by the Liouville theorem $Y(z)=\ID$ and the two Fuchsian systems 
must coincide. The proof in the case of $\mu=0$ is analogous. {\hfill QED}
\vskip 0.2 cm

\noindent{\bf Remark 5.}\quad The above argument fails for $R=0$, $\mu\neq0$.
Indeed, the fundamental matrices at $\infty$ of the form \eqrefp{M3}, with 
$R=0$, exist, but all the fundamental matrices of the form
$$
Y_\infty=\left(\ID +{\cal O}({1\over z})\right) z^{-{\cal A}_\infty} B
$$
with any constant matrix $B$, give the same monodromy matrix at infinity. 
Thus, chosen the fundamental matrices at $\infty$ of the two Fuchsian systems
as
$$
Y_\infty^{(i)}=\left(\ID +{\cal O}({1\over z})\right) z^{-{\cal A}_\infty}
B^{(i)}, \quad i=1,2,
$$
for some constant matrices $B^{(1)}$ and $B^{(2)}$, the above defined matrix 
$Y(z)$ is no more an analytic function near infinity and thus the uniqueness 
is not assured.
\vskip 0.2 cm

The theory of the isomonodromy deformations is described by the following two 
results (see [Sch]):

\proclaim Theorem 6. Let $M_1$, $M_2$, $M_3$ be the monodromy matrices of 
the Fuchsian system:
$$
\ddz Y^0=
\left({A^0_1\over z-u^0_1}+{A^0_2\over z-u^0_2}+
{A^0_3\over z-u^0_3}\right)Y^0,
\autoeqno{N7}
$$
of the above form \eqrefp{M2}, with $R\neq 0$, with pair-wise distinct poles, 
and with respect to some basis $\gamma_1,\gamma_2,\gamma_3$ of the loops in 
$\pi_1\left(\overline{\bf C}\backslash\{u^0_1,u^0_2,u^0_3,\infty\},
\infty \right)$. Then there exists a neighborhood $U\subset{\bf C}^3$ of 
the point $u^0=(u^0_1,u^0_2,u^0_3)$ such that, for any 
$u= (u_1,u_2,u_3)\in U$,
there exists a unique triple ${\cal A}_1(u)$, ${\cal A}_2(u)$, 
${\cal A}_3(u)$ of analytic matrix valued functions such that:
$$
{\cal A}_i(u^0)={\cal A}_i^0,\quad i=1,2,3,
$$
and the monodromy matrices of the Fuchsian system
$$
\ddz Y= A(z;u) Y=
\left({{\cal A}_1(u)\over z-u_1}+{{\cal A}_2(u)\over z-u_2}+
{{\cal A}_3(u)\over z-u_3}\right)Y, \autoeqno{N8}
$$
with respect to the same basis\footnote{${}^{1}$}{Observe that the 
basis $\gamma_1,\gamma_2,\gamma_3$ of 
$\pi_1\left(\overline{\bf C}\backslash\{u_1,u_2,u_3,\infty\},\infty \right)$ 
varies continuously with small variations of $u_1,u_2,u_3$. This new basis 
is homo-topic to the initial one, so we can identify them.}
 $\gamma_1,\gamma_2,\gamma_3$ of the loops, coincide with the given $M_1$, 
$M_2$, $M_3$.
The matrices ${\cal A}_i(u)$ are the solutions of the Cauchy problem with 
the initial data ${\cal A}_i^0$ for the following Schlesinger equations:
$$
{\partial\over\partial u_j}{\cal A}_i= 
{[{\cal A}_i,{\cal A}_j]\over u_i-u_j},\quad
{\partial\over\partial u_i}{\cal A}_i= 
-\sum_{j\neq i}{[{\cal A}_i,{\cal A}_j]\over u_i-u_j}. 
\autoeqno{N10}
$$
The solution $Y_\infty^0(z)$ of \eqrefp{N7}, with the form \eqrefp{M3},
can be, uniquely continued, for $z\not\in U$, to an analytic function 
$Y_\infty(z,u)$, $u\in U$, such that
$$
Y_\infty(z,u^0)=Y_\infty^0(z).
$$
Moreover the functions ${\cal A}_i(u)$ and $Y_\infty(z,u)$ can be continued 
analytically to global meromorphic functions on the universal coverings of
$$
{\bf C}^3\backslash\{diags\}:=
\left\{(u_1,u_2,u_3)\in{\bf C}^3\,|\,u_i\neq u_j\,\hbox{for}\, 
i\neq j\right\},
$$
and
$$
\left\{(z,u_1,u_2,u_3)\in{\bf C}^4\,|\,u_i\neq u_j\,\hbox{for}\, 
i\neq j\,\hbox{and}\,z\neq u_i,\, i=1,2,3\right\},
$$
respectively.

The proof of this theorem can be found, for example, in 
[Mal], [Miwa], [Sib]. 

\proclaim Theorem 7. Given three arbitrary non commuting matrices 
$M_1,M_2,M_3$ , satisfying \eqrefp{N5} and \eqrefp{N6}, with $M_\infty$ of 
the form \eqrefp{N9}, with $R\neq0$ and given a point
$u^0=(u^0_1,u^0_2,u^0_3)\in {\bf C}^3\backslash\{diags\}$, for any 
neighborhood $U$ of $u^0$, there exist $(u_1,u_2,u_3)\in U$ and a 
Fuchsian system of the form \eqrefp{M1}, with the given monodromy 
matrices, with the given $R$, with poles in $u_1,u_2,u_3$ and with a fixed 
value $\mu$ such that ${\rm Tr}M_\infty=2\cos\pi\mu$.

\noindent Proof. First, observe that if the matrices $M_1,M_2,M_3$, satisfying 
\eqrefp{N5} and \eqrefp{N6}, with $M_\infty$ of the form \eqrefp{N9} with 
$R\neq0$ commute then they are all lower triangular or all upper 
triangular. Since their eigenvalues are all equal to one, $M_\infty$ has 
eigenvalues equal to one, and thus $\mu\in\interi$.

Now, consider three arbitrary matrices $M_1,M_2,M_3$, 
satisfying \eqrefp{N5} and \eqrefp{N6}, with $M_\infty$ of the form 
\eqrefp{N9}, with $R\neq0$. In [Dek] it is proved that for any given 
point $u^0=(u^0_1,u^0_2,u^0_3)\in {\bf C}^3\backslash\{diags\}$, and for any 
neighborhood $U$ of $u^0$, there exist $(u_1,u_2,u_3)\in U$ and a Fuchsian 
system 
$$
\ddz Y=
\left({{\cal A}_1\over z-u_1} + {{\cal A}_2\over z-u_2} + 
{{\cal A}_3\over z-u_3} \right)  Y,
\qquad z\in\overline{\bf C}\backslash\{u_1,u_2,u_3,\infty\},
$$
with the given monodromy matrices, and with $\mu$ fixed up to $\mu\to\mu+n$, 
$n\in\interi$.

We want to build two gauge transformations which map the obtained Fuchsian 
system of the form \eqrefp{M1}, with some given non-commuting monodromy 
matrices and some value of $\mu$, to another Fuchsian system of the same form 
with the same monodromy matrices and with the value  $-\mu$ and $\mu+1$ 
respectively.

For $\mu\neq 0$ the constant gauge transformation
$$
G=\pmatrix{0&1\cr 1&0\cr}
$$
is such that the new Fuchsian system with 
$\tilde{\cal A}_i=G^{-1} {\cal A}_i G$,
has the same monodromy matrices $M_1,M_2,M_3$ and 
$$
\tilde{\cal A}_\infty=\pmatrix{-\mu&0\cr 0&\mu\cr}.
$$ 
So, the above gauge transformation maps the obtained Fuchsian system 
correspondent to the given monodromy matrices and some value of $\mu\neq0$ 
to another Fuchsian system of the same form with the same monodromy matrices 
and with the value  $-\mu$. Now, we want to build the analogous gauge 
transformation mapping $\mu$ to $\mu+1$.

First, observe that the matrices ${\cal A}_i$ can be parameterized as follows
$$
{\cal A}_i=\pmatrix{a_i b_i &-b_i^2\cr a_i^2 & -a_i b_i\cr},\autoeqno{p}
$$
for some $a_i,b_i\in\complessi$, $i=1,2,3$, with 
$$
\eqalign{
&\hbox{for}\quad\mu\neq0,\quad
\sum_{i=1}^3 a_i b_i=-\mu,\quad  \sum_{i=1}^3 a_i^2=\sum_{i=1}^3 b_i^2=0\cr
&\hbox{for}\quad\mu=0,\quad
\sum_{i=1}^3 a_i b_i=0,\quad  \sum_{i=1}^3 a_i^2=0,\quad
\sum_{i=1}^3 b_i^2=1.\cr}
$$
If $a_i=0$ (or $b_i=0$) for every $i=1,2,3$, then all the matrices 
${\cal A}_i$ are upper (resp. lower) triangular, then the matrices 
$M_1,M_2,M_3$ are upper (resp. lower) triangular, and thus commuting.
So, for a triple of non-commuting monodromy matrices at least one of the $a_i$
and one of the $b_i$ must be different from zero. Moreover, for $R\neq0$, and
for every $\mu$, $\sum a_i^2 u_i\neq 0$. In fact, if $\mu=-{1\over2}$, 
$\sum a_i^2 u_i=R_{21}\neq 0$. For $\mu\neq-{1\over2}$, if
$\sum a_i^2 u_i= 0$ then, being $a_1^2= -a_2^2-a_3^2$, one obtains 
$a_2^2=-a_3^2 {u_3-u_1\over u_2-u_1}$ and thus
$$
{\partial\over\partial u_1} a_2^2 = 
- {u_3-u_1\over u_2-u_1} {\partial\over\partial u_1} a_3^2
- a_3^2 {u_3-u_2\over(u_2-u_1)^2}.
$$
By the Schlesinger equations 
$$
2 a_1 b_1 a_2^2 -2 a_2 b_2 a_1^2 = 2 a_3 b_3 a_1^2 -2 a_1 b_1 a_2^2 
- a_3^2 {u_3-u_2\over u_2-u_1},
$$
and imposing $\sum a_i^2= 0$, $\sum a_i b_i= 0$, one obtains
$$
2 a_1^2\mu+a_3^2 {u_3-u_2\over u_2-u_1}=0
$$
that for $\mu=0$ leads to $a_3^2=0$ and thus $a_2=0$ and $a_1=0$, for 
$\mu\neq0$ leads to $a_1^2=-{a_3^2\over 2\mu}{u_3-u_2\over u_2-u_1}$. 
Imposing $\sum_{i=1}^3 a_i^2=0$, one obtains for $\mu\neq0$ 
$$
a_3^2 {(1+2\mu)(u_3-u_2)\over2\mu(u_2-u_1)}=0,
$$
that for $\mu\neq-{1\over2}$ implies $a_3^2=0$ and thus $a_2=0$ and $a_1=0$.
Analogously, one can show that for $\mu\neq0$ and $R\neq0$, 
$\sum b_i^2 u_i\neq 0$.

For $\mu\neq -1,0,-{1\over 2}$ the gauge transform $Y=G(z)\tilde  Y$ with
$$
G(z)=\pmatrix{1&0\cr 0& 0\cr} z+\pmatrix{{a\over b}&b\cr -{1\over b}&0\cr}
$$
for $b={2\mu+1\over \sum_{i=1}^3 a_i^2 u_i}$ and 
$a=-{b^2\over 2(1+\mu)}\left({2\over b}\sum_{i=1}^3 a_i b_i u_i+
\sum_{i=1}^3 a_i^2 u_i^2\right)$ is well defined because as observed above 
$\sum_{i=1}^3 a_i^2 u_i\neq0$ and it is such that the new Fuchsian system 
$$
\ddz\tilde Y=
\left({\tilde{\cal A}_1\over z-u_1} + {\tilde{\cal A}_2\over z-u_2} + 
{\tilde{\cal A}_3\over z-u_3} \right)\tilde  Y,
$$
with $\tilde{\cal A}_i=G(u_i)^{-1} {\cal A}_i G(u_i)$,
has the same monodromy matrices $M_1,M_2,M_3$ and 
$$
\tilde{\cal A}_\infty=\pmatrix{\mu+1&0\cr 0&-\mu-1\cr}.
$$ 
So, the above gauge transformation maps the obtained Fuchsian system 
corresponding to the given monodromy matrices and some value of 
$\mu\neq 0,-1,-{1\over 2}$ to another Fuchsian system of the same form with 
the same monodromy matrices and with the value $\mu+1$. 

In this way, all the half-integer values and all the non-zero integer values 
of the index $\mu$ are related via some gauge transformation. To conclude the 
proof, one has to consider the case of $\mu=0$. For a triple of non-commuting  
monodromy matrices with $\mu=0$, the gauge transformation $Y=G(z)\tilde Y$ with
$$
G(z)=\pmatrix{1&0\cr 0& 0\cr} z+\pmatrix{g_{11}&g_{12}\cr g_{21} &0\cr}
$$
with $g_{21}=-\sum a_i^2 u_i$, $g_{11}={1\over 2}\left(g_{21}-
2\sum_{i=0}^3 a_i b_i u_i+ {1\over g_{21}}\sum_{i=0}^3 a_i^2 u_i^2\right)$,
is well defined and it maps the Fuchsian system corresponding to the given 
triple of monodromy matrices to a new Fuchsian system with 
$\tilde{\cal A}_i=G(u_i)^{-1} {\cal A}_i G(u_i)$, with the same monodromy 
matrices $M_1,M_2,M_3$ and 
$$
\tilde{\cal A}_\infty=\pmatrix{1&0\cr 0&-1\cr}.
$$ 
In the same way, the gauge transformation $Y=G(z)\tilde Y$ with
$$
G(z)=\pmatrix{0&0\cr 0&1\cr} z+\pmatrix{0&g_{12}\cr g_{21} &g_{22}\cr}
$$
$g_{12}=-\sum b_i^2 u_i$, $g_{21}= \sum_{i=0}^3 a_i b_i u_i-
{1\over g_{12}}\sum_{i=0}^3 b_i^2 u_i^2$, and any $g_{22}\neq0$, is well 
defined and it maps any Fuchsian system with 
$$
\tilde{\cal A}_\infty=\pmatrix{1&0\cr 0&-1\cr},
$$ 
corresponding to the given triple of non-commuting monodromy matrices
to a new Fuchsian system with $\tilde{\cal A}_i=G(u_i)^{-1}{\cal A}_i G(u_i)$,
with the same monodromy matrices $M_1,M_2,M_3$ and 
$$
\tilde{\cal A}_\infty=\pmatrix{0&1\cr 0&0\cr}.
$$ 
This concludes the proof of the theorem. {\hfill QED}
\vskip 0.2 cm

\noindent{\bf Remark 6.} Existence statements of Theorems 6 and 7 can be 
proved also for triples of monodromy matrices such that $R=0$, but as 
stressed in Remark 5, uniqueness is lost.
\vskip 0.2 cm

Let me now explain, following [JMU], how to rewrite the Schlesinger equations 
\eqrefp{M2} in terms of the PVI$\mu$ equation. Observe that, for $\mu\neq0$,
the Schlesinger equations \eqrefp{M2} with fixed ${\cal A}_\infty$ are 
invariant with respect to the gauge transformations of the form:
$$
{\cal A}_i\mapsto D^{-1} {\cal A}_i D,\quad i=1,2,3,\quad\hbox{for any $D$ 
diagonal matrix}.\autoeqno{eqiuv}
$$ 
Such a diagonal conjugation changes the value of $R$. So, we introduce two 
coordinates $(p,q)$ on the quotient of the space of the matrices satisfying 
\eqrefp{N10} with respect to the equivalence relation \eqrefp{eqiuv}
and a coordinate $k$ that takes account of the changes of $R$ due to 
the above diagonal conjugations. The coordinate $q$ is the root of the 
following linear equation:
$$
[{\cal A}(q;u_1,u_2,u_3)]_{12}=0, 
$$
and $p$ and $k$ are given by:
$$
p=[{\cal A}(q;u_1,u_2,u_3)]_{11},\qquad
k=[{\cal A}(z;u_1,u_2,u_3)]_{12}{P(z)\over\mu(q-z)},
$$
where ${\cal A}(z;u_1,u_2,u_3)$ is given in \eqrefp{N8} and 
$P(z)=(z-u_1)(z-u_2)(z-u_3)$. The matrices ${\cal A}_i$ are uniquely 
determined by the coordinates $(p,q)$, and $k$ and expressed rationally in 
terms of them:
$$
\left\{
\eqalign{
\left( {\cal A}_i\right)_{11}&=-\left( {\cal A}_i\right)_{22}
={q-u_i\over2\mu P'(u_i)}\left[P(q)p^2+2\mu {P(q)\over q-u_i} p+
\mu^2(q+2u_i-\sum_ju_j)\right],\cr
\left( {\cal A}_i\right)_{12}&=-\mu k {q-u_i\over P'(u_i)},\cr
\left( {\cal A}_i\right)_{21}&= k^{-1} {q-u_i\over4\mu^3 P'(u_i)}
\left[P(q)p^2+2\mu {P(q)\over q-u_i} p
+\mu^2(q+2u_i-\sum_ju_j)\right]^2,\cr}\right.
\autoeqno{N12.5}
$$
for $i=1,2,3$, where $P'(z)={{\rm d}P\over{\rm d}z}$. 

The Schlesinger equations \eqrefp{N10} in the $(p,q,k)$ variables are
$$
\left\{\eqalign{
{\partial q\over\partial u_i} &=
{P(q)\over P'(u_i)}\left[2 p + {1\over q-u_i}\right]\cr
{\partial p\over\partial u_i} &=
-{P'(q) p^2 +(2q+u_i-\sum_ju_j)p+\mu(1-\mu)\over P'(u_i)},
\cr}\right.
\autoeqno{N13}
$$
and
$$
{\partial\log(k)\over\partial u_i}=(2\mu-1){q-u_i\over P'(u_i)}.
$$
for $i=1,2,3$. The system of the {\it reduced Schlesinger equations}\/ 
\eqrefp{N13} is invariant under the transformations of the form
$$
u_i\mapsto a u_i + b,\qquad
q\mapsto a q + b,\qquad
p\mapsto {p\over a}, \qquad\forall a,b\in\complessi,\quad a\neq0.
$$
We introduce the following new invariant variables:
$$
x={u_2-u_1\over u_3-u_1},\quad
y={q-u_1\over u_3-u_1};\autoeqno{N13.5}
$$
the system \eqrefp{N13}, expressed in the these new variables, reduces to 
the PVI$\mu$ equation for $y(x)$. 

The reduced Schlesinger equations admit the following {\it singular 
solutions}\/
$$
q=u_i\quad\hbox{for}\quad i=1,2,3.
$$
In [DM] is shown that $q=u_i$ if and only of the matrix ${\cal A}_i$ is 
identically equal to $0$, or, equivalently, the monodromy matrix $M_i$ is 
equal to the identity. The singular solutions do not give any solution of the 
PVI$\mu$ equation, while all the other solutions of the reduced Schlesinger 
equations do.

For $\mu=0$, the matrices ${\cal A}_i$ can be parameterized as in \eqrefp{p}
for some $a_i,b_i\in\complessi$. One can introduce the coordinates $(p,q)$ as 
above and obtain:
$$
a_1^2={p^2 P(q)\over\Delta}(u_3-u_2)(q-u_1),\quad b_1^2=0,\quad a_1 b_1=0,
$$
$$
a_2^2={p^2 P(q)\over\Delta}(u_1-u_3)(q-u_2),
\quad b_2^2={q-u_3\over q-u_1}{u_2-u_1\over u_2-u_3},
\quad a_2 b_2=p{(q-u_2)(q-u_3)\over u_2-u_3},
$$
$$
a_3^2={p^2 P(q)\over\Delta}(u_2-u_1)(q-u_3),
\quad b_3^2={q-u_2\over q-u_1}{u_3-u_1\over u_3-u_2},
\quad a_3 b_3=-p{(q-u_2)(q-u_3)\over u_2-u_3},
$$
where $P(q)=(q-u_1)(q-u_2)(q-u_3)$ and $\Delta=(u_3-u_2)(u_3-u_1)(u_2-u_1)$.
Introducing the variables $(y,x)$ as above, it is straightforward to verify 
that the Schlesinger equations for the matrices \eqrefp{p} are satisfied iff 
$y(x)$ satisfies PVI$_{\mu=0}$.

Observe that the Schlesinger equations for the matrices \eqrefp{p} admit the 
trivial solutions $a_i=0$, $\forall\,i=1,2,3$ or $b_i=0$, $\forall\,i=1,2,3$ 
correspondent respectively to the triple of commuting monodromy matrices
$$
M_i=\pmatrix{1 &-2\pi i b_i^2\cr 0&1},\quad\hbox{or}\quad
M_i=\pmatrix{1 &0\cr 2\pi i a_i^2 &1}.\autoeqno{M14}
$$
All the triples of commuting monodromy matrices can be realized by a trivial
solution to the Schlesinger equations for $\mu=0$. For such trivial solutions 
the coordinate $q$ is not defined. By the way, one cannot say that the 
triples of commuting monodromy matrices do not correspond to any solution of 
PVI$_{\mu=0}$. Indeed this equation coincides with PVI$_{\mu=1}$ and it is
not excluded that particular triples of commuting monodromy matrices are
realized by non-trivial Fuchsian system with $\mu=1$ (in Section 6, I show that
there exists a unique triple of this kind and determine it). 

\vskip 0.2 cm
Till now I supposed $R\neq 0$, because the uniqueness of the correspondence 
between triples of monodromy matrices and Fuchsian systems with a given set 
of poles is not assured for $R=0$. Thus I treat the case of $R=0$ separately.

\proclaim Lemma 11. For half integer values of $\mu$, the equation $R=0$ is 
satisfied if and only if the reduced Schlesinger equations give rise to Chazy 
solutions for any $\mu\neq{1\over 2}$ or to the singular solution
$q=\infty$ for $\mu={1\over 2}$.

\noindent Proof. Consider the case $\mu=-{1\over2}$ (as shown in Section 4,
all the other cases with half integer $\mu\neq{1\over2}$ are equivalent to 
it). The equation $R=0$ is satisfied iff 
$$
\left(A_1 u_1+A_2 u_2 +A_3 u_3 \right)_{21}=0.
$$
Writing the above equation in terms of $y$, $y'$ and $x$, one realizes that 
it coincides with the equation $Q(y,y',x)=0$ which is satisfied only by the
Chazy solutions.  


In the case $\mu ={1\over2}$, the equation $R=0$ is satisfied iff 
$$
\left(A_1 u_1+A_2 u_2 +A_3 u_3 \right)_{12}=0,
$$
that leads to the singular solution $q=\infty$. In fact, in the equation for $q$
$$
\left(A_1 (u_2+u_3)+A_2 (u_1+u_3)+A_3 (u_1+u_2) \right)_{12} q =
\left(A_1 u_2u_3+A_2 u_1u_3+A_3 u_1u_2\right)_{12},
$$
the coefficient of $q$ is zero and the right-hand side is non-zero because 
$(A_i)_{12}\neq 0$, $\forall\, i=1,2,3$. In fact if one of the $(A_i)_{12}$ 
is zero then, being $\sum(A_i)_{12}=0$, for $R=0$ all of them are $0$.
Requiring that the determinant of the matrices $A_i$ is zero, one obtains that
also the elements $(A_i)_11$ are zero, that is $\mu=0$ that leads a 
contradiction. {\hfill QED}

\proclaim Lemma 12. The equation $R=0$ is not satisfied on any solution
of the reduced Schlesinger equations for integer $\mu$.

\noindent Proof. Consider the cases of integer $\mu$. They can all be treated 
as the case $\mu=1$. In fact the case $\mu=0$ gives rise to the same PVI 
equation as the case $\mu=1$. Moreover all the other integer
values of $\mu$ are related to $\mu=1$ via birational canonical 
transformations of the form \eqrefp{S1}, the denominator of which never 
vanishes. For $\mu=1$ the equation $R=0$ in the $(p,q)$ coordinates is
$$
R_{12}=k\left[p^2 (q-u_1)(q-u_2)(q-u_3)-2 q+2 u_1+2 u_2+2 u_3\right],
\autoeqno{abs}
$$
which is never zero on the solutions of the reduced Schlesinger equations.
In fact from \eqrefp{abs}, one obtains:
$$
p^2={2(q-u_1-u_2-u_3)\over(q-u1)(q-u2)(q-u3)}.
$$
Differentiating both sides with respect to $u_i$ for all $i=1,2,3$, and 
substituting the reduced Schlesinger equation for 
${\partial p\over\partial u_i}$, one obtains
$$
p=-{q-u_1-u_2\over(q-u_1)(q-u_2)}=-{q-u_1-u_3\over(q-u_1)(q-u_3)}=
-{q-u_3-u_2\over(q-u_3)(q-u_2)},
$$
which can be satisfied only for $u_1=u_2=u_3$. 
{\hfill QED}

\vskip 0.2 cm
As shown above, the case $R=0$ can be realized only for half integer $\mu$.
It gives rise to the singular solution $q=\infty$ in the case of 
$\mu={1\over2}$ and to the Chazy solutions in the case of half-integer 
$\mu\neq{1\over2}$. 

\noindent We resume the results of the section in the following:

\proclaim Theorem 8. Given any triple of non-commuting monodromy matrices 
$M_1$, $M_2$, $M_3$ satisfying \eqrefp{N5} and \eqrefp{N6}  with 
$M_\infty$ given by \eqrefp{N9}, for $R\neq0$, none of them being 
equal to $\ID$, considered modulo diagonal conjugations, there exists unique 
branch of a non-Chazy solution to the PVI$\mu$ equation near a given point 
$x_0\in\overline{\bf C}\backslash\{0,1,\infty\}$ which defines a Fuchsian 
system of the form \eqrefp{M1} with the prescribed monodromy matrices $M_1$, 
$M_2$, $M_3$. Vice versa, given any branch of a non-Chazy solution to the 
PVI$\mu$ equation near a given point 
$x_0\in\overline{\bf C}\backslash\{0,1,\infty\}$, the correspondent triple 
of monodromy matrices $M_1$, $M_2$, $M_3$ satisfying \eqrefp{N5} and 
\eqrefp{N6}  with $M_\infty$ given by \eqrefp{N9}, for $R\neq0$, none of them 
being equal to $\ID$ is unique modulo diagonal conjugations.

Observe that permutations of the poles $u_i$ induce transformations of $(y,x)$ 
of the type $x\to1-x$, $y\to1-y$ and $x\to{1\over x}$ and $y\to{y\over x}$ 
and their compositions. These transformations preserve the PVI$\mu$ equation.

\vfill \eject

\semiautosez{6}{6. The structure of the analytic continuation.}

We parameterized branches of the non-Chazy solutions of PVI$\mu$ by triples of 
monodromy matrices. Recall that, according to Theorem 6, the 
solutions of PVI$\mu$, defined in a neighborhood of a given point 
$x_0\in\overline{\bf C}\backslash\{0,1,\infty\}$, can be analytically 
continued to a meromorphic function on the universal covering of 
$\overline{\bf C}\backslash\{0,1,\infty\}$. In [DM] it is shown that 
the procedure of the analytic continuation is described by the action of 
the pure braid group with three strings (see [Bir]),
$P_3=\pi_1\left({\bf C}^3\backslash\{diags\};\, u^0\right)$, on  
$\pi_1\left(\overline{\bf C}\backslash\{0,1,\infty\}\right)$. To simplify 
the computations, the procedure of the analytic continuation is extended 
to the full braid group $B_3$ that admits a presentation with generators 
$\beta_1$ and $\beta_2$ shown in figure 2, and defining relations
$\beta_1 \beta_2 \beta_1=\beta_2 \beta_1\beta_2$.

\midinsert
\centerline{\psfig{file=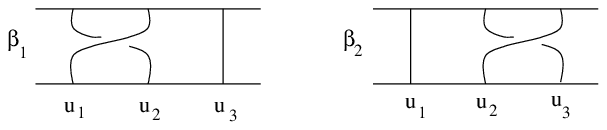,height=2cm}}\vskip 0.5 cm
\centerline{{\bf Fig.5.} The geometric representation of the generators of 
the braid group $B_3$.}
\endinsert

The action of $B_3$ on the generating loops of 
$\pi_1\left(\overline{\bf C}\backslash\{0,1,\infty\}\right)$ can be 
expressed in terms of monodromy matrices:

\proclaim Lemma 13. For the generators $\beta_1$, $\beta_2$ shown in the 
figure 2, the matrices $M_i^\beta$ have the following form:
$$
M_1^{\beta_1}=M_2,\quad M_2^{\beta_1}=M_2 M_1 M_2^{-1},
\quad M_3^{\beta_1}= M_3,\autoeqno{D2}
$$
$$
M_1^{\beta_2}=M_1,\quad M_2^{\beta_2}=M_3,
\quad M_3^{\beta_2}= M_3 M_2 M_3^{-1}.\autoeqno{D3}
$$

\noindent The proof can be found in [DM].
\vskip 0.2 cm

The action \eqrefp{D2}, \eqrefp{D3} of the braid group on the triples 
of monodromy matrices commutes with the diagonal conjugation of them; 
moreover the class of the singular solutions is closed under the analytic 
continuation. In fact if some of the matrices $M_i$ is equal to $\ID$, 
then for any $\beta$ there is a $j$ such that $M_j^\beta=\ID$. Moreover $R$ 
is preserved and the class of non-Chazy solutions is closed under the
analytic 
continuation. As a consequence the structure of the analytic continuation of 
the non-Chazy solutions of the PVI$\mu$ equation is determined by the action 
\eqrefp{D2}, \eqrefp{D3} of the braid group on the triples of monodromy 
matrices. 

I want to introduce a parameterization of the monodromy matrices and
write this action in terms of the parameters in the space of the monodromy 
data. I follow the same procedure of [DM], suitably modified due to the 
resonant value of $\mu$. The following lemmas can be proved as in [DM] 
(see lemma 1.4 and 1.5). 

\proclaim Lemma 14. Let ${\cal M}_1$, ${\cal M}_2$ and ${\cal M}_3$ be 
three linear non commuting operators 
${\cal M}_i:{\bf C}^2\rightarrow{\bf C}^2$ satisfying \eqrefp{N5} and 
\eqrefp{N6}, with ${\cal M}_\infty$ given by \eqrefp{N9} with $R\neq 0$. If 
two of the following numbers
$$
{\rm Tr}({\cal M}_1{\cal M}_2),\quad {\rm Tr}({\cal M}_1{\cal M}_3),
\quad {\rm Tr}({\cal M}_3{\cal M}_2)\autoeqno{M15}
$$
are equal to $2$, then one of the matrices of $M_i$ is equal to one.

\proclaim Lemma 15. Let ${\cal M}_1$, ${\cal M}_2$, ${\cal M}_3$ as in 
lemma 14.
\item{i)} If ${\rm Tr}({\cal M}_1{\cal M}_2)\neq 2$, then  there exists a 
basis in $\complessi^2$ such that, in this basis, the matrices $M_1$, 
$M_2$ and $M_3$ have the form
$$
M_1=\pmatrix{1&-x_1\cr 0&1\cr},\quad M_2=\pmatrix{1&0\cr x_1&1\cr},
\quad M_3=\pmatrix{1+ {x_2 x_3\over x_1}& -{x_2^2\over x_1}\cr 
{x_3^2\over x_1}&1-{x_2 x_3\over x_1}\cr},
\autoeqno{N15}
$$
where
$$
{\rm Tr}(M_1M_2)=2-x_1^2,\quad {\rm Tr}(M_3M_2)=2-x_2^2,\quad 
{\rm Tr}(M_1M_3)=2-x_3^2,
$$
and
$$
x_1^2+x_2^2+x_3^2-x_1 x_2 x_3=4\sin^2\pi\mu.\autoeqno{N16}
$$
\item{ii)} If two triples of matrices $M_1$, $M_2$, $M_3$ and $M'_1$, 
$M'_2$, $M'_3$, satisfying \eqrefp{N6}, with none of them equal to $\ID$, 
have the form \eqrefp{N15} with parameters $(x_1,x_2,x_3)$ and 
$(x'_1,x'_2,x'_3)$ respectively, then these triples are conjugated 
$$
M_i=T^{-1} M_i' T
$$
with some invertible matrix $T$, if and only if the triple 
$(x'_1,x'_2,x'_3)$ is equal to the triple $(x_1,x_2,x_3)$, up to the 
change of the sign of two of the coordinates.

Let me stress that for $\mu\in\interi$ the monodromy matrices can commute,
i.e. they can be of the form \eqrefp{M14}. The action \eqrefp{D2}, \eqrefp{D3} 
of the braid group simply permutes them. Thus the action \eqrefp{D2}, 
\eqrefp{D3} of the braid group does not mix the commuting triples \eqrefp{M14}
with the ones admitting the parameterization \eqrefp{M15}.

\proclaim Lemma 16. There exists only a one-parameter family of triples of
commuting monodromy matrices which give rise to solutions to PVI$\mu$ equation.
They are (up to diagonal conjugations)
$$
M_1=\pmatrix{1& i\pi a\cr 0 &1\cr},\quad
M_2=\pmatrix{1& i\pi (1-a)\cr 0 &1\cr},\quad
M_3=\pmatrix{1& i\pi \cr 0 &1\cr},\autoeqno{M20}
$$
where $a$ is a parameter. The correspondent solutions to PVI$_{\mu=1}$ 
consist of a one parameter family of rational solutions of the form: 
$$
y(x) = {a x\over1-(1-a)x},\qquad\hbox{for}\quad a\neq 0.\autoeqno{M21}
$$

\noindent Proof. Consider a triple of commuting monodromy matrices. As shown 
above, they are necessarily of the form \eqrefp{M14}, i.e. either all upper 
triangular or all lower triangular. Then the corresponding Fuchsian system 
admits a single-valued solution $Y(z)$. For $\mu=1$ (as shown 
in the proof of Lemma 12, all the other cases with integer $\mu$ are 
equivalent to this) such solution has only a pole of order one at infinity, 
i.e.
$$
Y(z)=\pmatrix{a z+b\cr c z+d \cr}.
$$
for some $a,b,c,d\in\complessi$. Substituting $Y$ in the Fuchsian system, 
one obtains $a=0$, $b=c k$, $d={c\over 2}(q-u_1-u_2-u_3)$ and $c\neq0$ iff 
$p\equiv0$. By direct substitution in the reduced Schlesinger equations, 
one can compute $q$ and determine the explicit form of Fuchsian system. Thus 
it is straightforward to compute the monodromy matrices, which turn out to 
have the form \eqrefp{M20}. Their orbit under the action of the braid group 
\eqrefp{D2}, \eqrefp{D3} consists of one point, up to permutations. Thanks to 
Theorem 8, the correspondent solution of PVI$\mu$, $\mu\in\interi$ consists 
only of one branch, i.e. it is rational and it is easy to see that, being 
$p\equiv 0$ it has the form \eqrefp{M21}.
{\hfill QED}
\vskip 0.2 cm

For the case when at most one of the numbers \eqrefp{M15} is equal to $2$,
or equivalently the monodromy matrices do not commute, one gives the following:

\vskip 0.2 cm
\noindent{\bf Definition:}
a triple $(x_1,x_2,x_3)$ is called {\it admissible}\/ if it has at most 
one coordinate equal to zero. Two triples are called 
{\it equivalent}\/ if they are equal up to the change of two signs 
of the coordinates. The equivalence class of a triple $(x_1,x_2,x_3)$ is
denoted by $[x_1,x_2,x_3]$
\vskip 0.2 cm

Observe that for an admissible triple $(x_1,x_2,x_3)$, none of the 
matrices \eqrefp{N15} is equal to the identity. So the admissible 
triples give rise to non-singular solutions of the reduced Schlesinger 
equations. Moreover two equivalent triples such that $R\neq0$ generate the 
same solution. 
 
\proclaim Theorem 9. In the case of $\mu\in\interi$ there exists a one 
parameter family of rational solutions of the form \eqrefp{M21}. The 
one-parameter family of Chazy solutions in the case of half-integer $\mu$, 
$\mu\neq{1\over2}$, corresponds to $[2,2,2]$. All the other solutions of 
PVI${\mu}$, with $2\mu\in\interi$, have branches which, near a 
given point $x_0\in\overline\complessi\backslash\{0,1,\infty\}$, are in 
one-to-one correspondence with the equivalence classes of the admissible 
triples satisfying \eqrefp{N16}. The structure of the analytic continuation of
them is determined by the action 
$$
\eqalign{
\beta_1:(x_1,x_2,x_3)&\mapsto(-x_1,x_3-x_1x_2,x_2),\cr
\beta_2:(x_1,x_2,x_3)&\mapsto(x_3,-x_2,x_1-x_2 x_3).\cr}\autoeqno{BD}
$$
of the braid group on the triples $(x_1,x_2,x_3)$.

\noindent Proof. The first claim is proved in Lemma 16.
The second claim follows 
from the fact that as proved in Lemma 11, Chazy solutions correspond to the 
case $R=0$, i.e. to $M_\infty=-\ID$. Since $M_\infty$ is invariant with 
respect to conjugations, it must be equal to $-\ID$ also in the canonical form 
\eqrefp{N15}. Solving the equations in $(x_1,x_2,x_3)$, one obtains that the 
triple of monodromy matrices is given by
$$
M_1=\pmatrix{1&-2\cr 0&1\cr},\quad
M_2=\pmatrix{1&0\cr 2&1\cr},\quad
M_3=\pmatrix{3&-2\cr 2&-1\cr}\autoeqno{BD1}
$$
The third and fourth claim of this theorem were proved in [DM]. {\hfill QED}
\vskip 0.2 cm

\semiautosez{7}{7. Solutions of the PVI$\mu$ equation having finite branching
and their monodromy data.}

In this section I classify all the monodromy data corresponding to solutions 
with a finite number of branches. The solutions corresponding to commuting
triples of monodromy matrices were already found in Lemma 16. All the 
other triples admit the parameterization \eqrefp{N15}, and the strategy is 
essentially the same of [DM]. I recall here the ideas and results of [DM] 
omitting all the proofs. 

Let $y(x)$ be a finite-branching solution. According to the 
Painlev\'e property, the ramification points of $y(x)$ are allowed to lie 
only at $0,\,1,\,\infty$ and the correspondent monodromy matrices, defined 
modulo diagonal conjugations, have a finite orbit under the action of the 
braid group \eqrefp{D2}, \eqrefp{D3}.

\proclaim Lemma 17. An admissible triple $(x_1,x_2,x_3)\not\in[2,2,2]$, 
specifies a finite-branching solution of PVI$\mu$, for $\mu$ given by 
\eqrefp{N16}, if and only if its orbit, under the action \eqrefp{BD} of the 
braid group, is finite.

\noindent{\bf Remark 7.}\quad 
Observe that evenif the orbit of the triple $(2,2,2)$ consists only of one 
point, up to equivalence, it gives rise to a one parameter family of 
transcendental solutions (the Chazy solutions). This is not surprising 
because this triple corresponds to the case $R=0$. As stressed in Lemma 10, 
the uniqueness of the Fuchsian system associated to triples of monodromy 
matrices in not assured in this case. Observe that the first term of 
the asymptotic behaviour of the Chazy solutions does not depend on the 
parameter $\nu$. This fact could be related to the fact that the monodromy 
matrices associated to the Chazy solutions do not depend on $\nu$.
\vskip 0.2 cm

Due to Lemma 17, the problem of the classification of the regular 
solutions of the PVI$\mu$ reduces again to the problem of the classification 
of all the finite orbits of the action \eqrefp{BD} under the braid group in 
the three dimensional space.  The following simple necessary condition for 
a triple $(x_1,x_2,x_3)$ to belong to a finite orbit is proved in [DM].

\proclaim Lemma 18. Let $(x_1,x_2,x_3)$ be a triple belonging to a 
finite orbit. Then:
$$
x_i=-2\cos\pi r_i,\quad r_i\in{\bf Q},\quad 0\leq r_i\leq 1.
\quad i=1,2,3,
\autoeqno{N17}
$$
here ${\bf Q}$ is the set of the rational numbers.

\noindent{\bf Remark 8.}\quad Thanks to the above lemma, for the finite orbits 
of the braid group, it is equivalent to deal with the triples 
$(x_1,x_2,x_3)$, or with the {\it triangles}\/ of angles 
$(\pi r_1,\pi r_2,\pi r_3)$, with $x_i=-2\cos\pi r_i$ and $0\leq r_i\leq 1$,
$r_i\in{\bf Q}$.

\vskip 0.3 cm
\noindent{\bf 7.1. Classification of the monodromy data corresponding 
to finite-branching solutions: case of $\mu\in\interi$.}\quad In this case 
the classification theorem of [DM] (see Theorem 1.6) is still valid, so there 
are no triples $(x_1,x_2,x_3)$ having a finite orbit. As a consequence, the 
only finite-branching solutions are the ones corresponding to
triples of monodromy matrices admitting the parameterization \eqrefp{M14}. 

\proclaim Lemma 19. The only finite-branching solutions of PVI$\mu$, for 
$\mu=1$, consist of a one-parameter family of rational solutions of the form 
\eqrefp{M21}. 

All the other cases of PVI$\mu$, for $\mu\in\interi\backslash\{0\}$ can be 
obtained from it via birational canonical transformations \eqrefp{S1}, the 
denominator of which never vanishes. The case $\mu=0$ is the same as $\mu=1$.

\vskip 0.3 cm
\noindent{\bf 7.2. Classification of the triples $(x_1,x_2,x_3)$ corresponding 
to finite-branching solutions: case of half-integer $\mu$.}\quad
For $\mu+{1\over2}\in\interi$, the triples $(x_1,x_2,x_3)$ satisfy
$$
x_1^2+x_2^2+x_3^2-x_1 x_2 x_3 =4,
$$ 
and the correspondent triangle is flat, $r_1+r_2+r_3=1$.

\proclaim Lemma 20. For every admissible triple $(x_1,x_2,x_3)$, 
$x_i=-2\cos\pi r_i$, with $r_i\in{\bf Q}$, the correspondent solution
of PVI$\mu$, for $\mu$ half-integer, is finite-branching.

\noindent Proof.
The action of the braid group on flat triangles can be written in the form
$$
\eqalign{
\beta_1:(r_1,r_2,r_3)&\mapsto (|1-r_1|,|r_1-r_2|, r_2),\cr
\beta_2:(r_1,r_2,r_3)&\mapsto (r_3,|1-r_2|,|r_3- r_2|).\cr}
\autoeqno{M16}
$$
As a consequence, it maps triangles with rational angles to triangles with 
rational angles. Moreover all the orbits are finite.
In fact, let $r_i={p_i\over q_i}$, for $p_i,q_i\in\interi$, $p_i<q_i$, 
$i=1,2,3$ and $n$ be the smallest common factor of $q_1,q_2,q_3$. The 
action of the braid group \eqrefp{M16} does not increase $n$, and all
the images of $(r_1,r_2,r_3)$ have the form $\left({\tilde p_1\over n}, 
{\tilde p_2\over n}, {\tilde p_3\over n}\right)$, with $\tilde p_i<n$. 
The number of possible triples of this kind is trivially finite.
{\hfill QED}

\vskip 0.2 cm
Recall that since the Chazy solutions are transcendental, the finite-branching 
solutions in the case of half-integer $\mu$, are necessarily of Picard type,
and thus they have asymptotic behaviour of algebraic type, see \eqrefp{P1.5}.
In the following Proposition, I relate the parameters $\nu_{1,2}$ of Picard 
solutions to the triples $(r_1,r_2,r_3)$.

\proclaim Proposition 2. i) The monodromy matrices corresponding to a solution 
$y(x)$ of the form \eqrefp{P1} are given by 
$$
M_1=\pmatrix{1&-2\cos\left({\pi\nu_2\over2}\right)\cr 0&1\cr},
\quad
M_2=\pmatrix{1&0\cr 2\cos\left({\pi\nu_2\over2}\right)&1\cr},
$$
$$
M_3=\pmatrix{1+{2\cos\left({\pi\nu_1\over2}\right)
\cos\left({\pi(\nu_1-\nu_2)\over2}\right)
\over\cos\left({\pi\nu_2\over2}\right)}
&{-2 \cos^2\left({\pi\nu_1\over2}\right)
\over\cos\left({\pi\nu_2\over2}\right)}\cr
{2 \cos^2\left({\pi(\nu_1-\nu_2)\over2}\right)
\over\cos\left({\pi\nu_2\over2}\right)}&
1-{2\cos\left({\pi\nu_1\over2}\right)
\cos\left({\pi(\nu_1-\nu_2)\over2}\right)
\over\cos\left({\pi\nu_2\over2}\right)}\cr}.
$$
That is the parameters $\nu_{1,2}$ and the triples $(x_1,x_2,x_3)$ are 
related as follows:
$$
\eqalign{
r_1 = {\nu_2\over2},\quad
r_2 = 1-{\nu_1\over2},\quad
r_3 = {\nu_1-\nu_2\over2}, &\qquad\hbox{for}\quad \nu_1>\nu_2,\cr
r_1 = 1-{\nu_2\over2},\quad
r_2 = {\nu_1\over2},\quad
r_3 = {\nu_2-\nu_1\over2},&\qquad\hbox{for}\quad \nu_1<\nu_2,\cr}
\autoeqno{M19}
$$
and viceversa
$$
\nu_1 = 2-2 r_2,\quad
\nu_2 = 2 r_1,\autoeqno{M19.5}
$$
where $x_i=-2\cos\pi r_i$.
ii) All the finite-branching solutions are indeed algebraic.
iii) The action of the braid group $B_3$ (pure braid group $P_3$) on 
$(x_1,x_2,x_3)$ corresponds to the action of $\Gamma$ ($\Gamma(2)$) on 
$(\nu_1,\nu_2)$.

\noindent Proof. i) The relation between the parameters $\nu_{1,2}$ and the 
exponents $(l_0,l_1,l_\infty)$ of the asymptotic behaviour was derived in 
Lemma 2.
It remains to find the relation between the exponents and the triangles. 
For $l_i\neq0$ for any $i=0,1,\infty$, this was already done in [DM] for 
PVI$\mu$, $2\mu\not\in\interi$. This relation is
$$
l_i=\left\{\eqalign{
2r_i\quad\hbox{for}\quad 0<r_i\leq{1\over 2},\cr
2-2r_i\quad\hbox{for}\quad {1\over 2}\geq r_i<1.\cr}
\right.\autoeqno{M17}
$$
This result can be extended to the case of $\mu={1\over2}$. Indeed, the 
procedure of reduction of the Fuchsian system to the systems $\hat\Sigma$ 
and $\tilde\Sigma$ of [DM] is the same. They are again reduced to the Gauss 
equation. 
The only difference appears in the computation of the connection matrices 
of the system $\hat\Sigma$. In fact, for $\mu={1\over2}$, the fundamental 
matrix at infinity is a Jordan block and has logarithmic type behaviour. 
The computations of the analytic continuation can be performed 
following the formulae of [Nor] and the connection 
matrices are computed as in [DM]. Then, using \eqrefp{M17}, \eqrefp{P2.5} and 
\eqrefp{P2.55}, one can show \eqrefp{M19}, for $\nu_i\neq0$.
Let us suppose $\nu_i=0$ for some $i$, for example $\nu_2=0$, i.e. $l_0=0$.
Then $\nu_1\neq0$, and $l_1=l_\infty$, $l_0=0$. I can take 
$r_2=1-r_3={l_1\over2}$. Since $r_1+r_2+r_3=1$, $r_1$ must be $0$, and the 
lemma is proved also for $\nu_2=0$. 

ii) Follows from the fact that finite-braching solutions correspond to rational
values of $(r_1,r_2,r_3)$, i.e. to rational $(\nu_1,\nu_2)$, thus they are 
algebraic (see Lemma 3).

iii) The fact that the action of the braid group $B_3$ (pure braid group $P_3$)
on $(x_1,x_2,x_3)$ corresponds to the action of $\Gamma$ ($\Gamma(2)$) on 
$(\nu_1,\nu_2)$ is easily derived by the formulae \eqrefp{M19}, \eqrefp{M19.5}
relating $(\nu_1,\nu_2)$ and the angles, and by the formula \eqrefp{M16}.
{\hfill QED}

\vskip 0.2 cm
\noindent{\bf Remark 9.}\quad Observe that in the limit $\nu_{1,2}\to0$ 
the above formulae for the monodromy matrices give the matrices \eqrefp{BD1}. 
Indeed, as shown in Section 4.1, Chazy solutions are limits of Picard type 
solutions for $\nu_{1,2}\to0$. 

\vskip 0.3 cm
\semiautosez{8}{8. Algebraic solutions and finite irreducible reflection 
groups.}

We reformulate here the above parameterization of the monodromy data, for 
the case of half-integer $\mu$, by flat triangles $(r_1,r_2,r_3)$, or 
equivalently,
by couples of constants $(\nu_1,\nu_2)$ in a more geometric way.
Consider the Euclidean three-dimensional space $E$ and three planes
$(p_1,p_2,p_3)$ all intersecting in one point, let $(r_1,r_2,r_3)$ be
the angles between them and $(e_1,e_2,e_3)$ the vectors normal to them.
Define three reflections $R_1,R_2,R_3$ with respect to the three planes 
$(p_1,p_2,p_3)$:
$$
R_i:\eqalign{
E&\rightarrow E\cr
x&\mapsto x-(e_i,x)e_i\cr}
\qquad i=1,2,3.
$$
Let us consider the group $G$ of the linear transformations of $E$, generated 
by the three reflections $R_1,R_2,R_3$. Its Gram matrix is
$$
g:=\pmatrix{2&x_1&x_3\cr
x_1&2&x_2\cr
x_3&x_2&2\cr}\autoeqno{N29}
$$
where $x_i=-2\cos\pi r_i$. Observe that $g$ is always singular:
$$
\det g = 8- 2(x_1^2+x_2^2+x_3^2-x_1 x_2 x_3)= 8\cos^2\pi \mu=0,
$$
then the normal vectors $(e_1,e_2,e_3)$ are linearly dependent.

Observe that the group $G$ is invariant under analytic continuation. In fact, 
as shown in [DM], the action of the braid group on the triples $(x_1,x_2,x_3)$
can be interpreted as an action on the correspondent generating reflections
$$
\eqalign{
\beta_1:(R_1,R_2,R_3)\mapsto 
(R_1,R_2,R_3)^{\beta_1}:=&(R_2,R_2R_1R_2,R_3),\cr
\beta_2:(R_1,R_2,R_3)\mapsto 
(R_1,R_2,R_3)^{\beta_2}:=&(R_1,R_3, R_3R_2R_3),\cr}\autoeqno{B19}
$$
where $\beta_{1,2}$ are the standard generators of the braid group. The 
groups generated by the reflections $(R_1,R_2,R_3)$ and 
$(R_1,R_2,R_3)^{\beta}$ coincide for any $\beta\in B_3$. In particular the 
following lemma holds true:

\proclaim Lemma 21. For any braid $\beta\in B_3$, the images 
$\beta(R_1,R_2,R_3)$ are reflections with respect to some planes orthogonal 
to some new basic vectors $(e^\beta_1,e^\beta_2,e^\beta_3)$. The Gram matrix 
with respect to the basis $(e^\beta_1,e^\beta_2,e^\beta_3)$ has the form:
$$
(e^\beta_i,e^\beta_i)=2,\quad i=1,2,3,\quad
(e^\beta_1,e^\beta_2)=x^\beta_1,\quad (e^\beta_2,e^\beta_3)=x^\beta_2,
\quad (e^\beta_1,e^\beta_3)=x^\beta_3,
$$
where $(x^\beta_1,x^\beta_2,x^\beta_3)=\beta(x_1,x_2,x_3)$.

Consider any algebraic solution of PVI$\mu$ with half-integer $\mu$. 
According to Lemma 3, it is specified by a couple of coprime integers 
$0\leq M<N$, and, thanks to Proposition 2, the correspondent triangles belong 
to the orbit of $\left(0,{M\over 2N},1-{M\over 2N}\right)$. Thus two mirrors 
coincide and form an angle ${\pi M\over 2 N}$ with the third. This is called 
a {\it dihedral kaleidoscope}\/ (sse [Cox]). The generated group is the 
dihedral group $D(\hat N)$ realized as symmetry group of a regular 
star-polygon with $\hat N$ edges and density $\hat M$, where
$$
\hat N=\left\{\eqalign{N\quad\hbox{if}\quad & M\quad\hbox{is even},\cr
2N\quad\hbox{if}\quad & M\quad\hbox{is odd},\cr}\right.\qquad
\hat M=\left\{\eqalign{{M\over2}\quad\hbox{if}\quad & M\quad\hbox{is even},\cr
M\quad\hbox{if}\quad & M\quad\hbox{is odd}.\cr}\right.
$$
Resuming, I proved the following

\proclaim Theorem 10. The algebraic solutions of PVI$\mu$ with any 
half-integer $\mu$, are in one to one correspondence with regular polygons 
and star-polygons in the plane.

\vskip 0.2 cm
\noindent{\bf Remark 10.}\quad The algebraic solutions \eqrefp{a2}, 
\eqrefp{b2} and \eqrefp{g2} correspondent to the values $N=3$ and $M=2$, 
$N=2$ and $M=1$, $N=3$ and $M=1$ respectively, correspond to $D(3)$, $D(4)$, 
$D(6)$ which coincide with $A_2$, $B_2$ and $G_2$ respectively 
\vskip 0.2 cm

In the Chazy case, all the mirrors coincide and the generated group is the 
cyclic group of order two, the symmetry group of a dighon.

For integer $\mu$ the only algebraic solutions are the ones correspondent to
the triple $(0,0,0)$. The correspondent Gram matrix is 
$$
g=\pmatrix{2&0&0\cr 0&2&0\cr 0&0&2\cr}.
$$
This means that the basic vectors are all orthoganal one to eachother, i.e. 
they form a Cartesian frame, the correspondent reflection group is generated 
by the three inversions of the Cartesian coordinates and it is abelian.

\vskip 0.3 cm
\baselineskip=12pt
\noindent{\bf BIBLIOGRAPHY }
\vskip 0.3 cm

\item{[Bir]}
J.S. Birman, {\it Braids, Links, and Mapping Class groups,}\/ Ann. Math. 
Stud. Princeton University (1975).
\vskip 0.2 cm

\item{[Cha]}
J. Chazy, Sur les equations differentielles dont l'integrale generale
possede un coupure essentielle mobile, {\it  C.R. Acad. Paris}\/ {\bf 150} 
(1910), 456-458.
\vskip 0.2 cm

\item{[Dek]}
W. Dekkers, {\it The matrix of a connection having regular singularities 
on a vector bundle of rank $2$ on $P^1({\bf C})$,}\/ Lec. Notes Math. 
{\bf 712}, p. 33-43 (1979).
\vskip 0.2 cm

\item{[Dub]} 
B. Dubrovin, {\it Geometry of 2D Topological Field Theories,}\/ Lect. 
Notes Math. {\bf 1620}, (1996).
\vskip 0.2 cm

\item{[DM]} 
B. Dubrovin and M. Mazzocco, {\it Monodromy of certain Painlev\'e 
transcendents and reflection groups,}\/ SISSA preprint n.149/97/FM 
(1997).
\vskip 0.2 cm

\item{[ItN]} 
A.R. Its and V.Yu. Novokshenov, {\it The isomonodromic deformation method 
in the theory of Painlev\'e equations,}\/ Lec. Not. Math. {\bf 1191} 
(1986) Springer.
\vskip 0.2 cm

\item{[FlN]}
H. Flashka and A.C. Newell, Monodromy and spectrum preserving deformations,
{\it Comm. Math. Phys.}\/ {\bf 76} (1980).
\vskip 0.2 cm

\item{[Fuchs]}
R. Fuchs, \"Uber lineare homogene Differentialgleichungen zweiter 
Ordnung mit im drei im Endrichen gelegene wesentlich singul\"aren 
Stellen, {\it Math. Ann.}\/ {\bf 63} (1907) 301-321.
\vskip 0.2 cm

\item{[Gamb]}   
B. Gambier, Sur les Equations Differentielles du Second Ordre et du 
Primier Degr\`e dont l'Integrale est a Points Critiques Fixes, {\it Acta 
Math.}\/  {\bf 33}, (1910) 1-55.
\vskip 0.2 cm

\item{[Hal]}
G.H. Halphen, Sur un systeme d'equations differentielles, {\it C.R. Acad
Sc. Paris}\/ {\bf 92} (1881) 1001-1007.
\vskip 0.2 cm

\item{[Hit]}
N.J. Hitchin, Twistor spaces, Einstein metrics and isomonodromic 
deformations, {\it J. Differential Geometry}\/ {\bf 42}, No.1 July 1995.
\vskip 0.2 cm

\item{[Ince]}
E.L. Ince, {\it Ordinary Differential Equations,}\/
Dover Publications, New York (1956).
\vskip 0.2 cm

\item{[Man]}
Yu.I. Manin, {\it Sixth Painlev\'e equation, universal elliptic curve and 
mirror of ${\bf P}^2$,}\/ preprint (1996).
\vskip 0.2 cm

\item{[Mal]}
B. Malgrange, Sur les deformations isomonodromiques I, singularit\'es 
r\'eguli\`eres, {\it Seminaire de l'Ecole Normale Superieure 1979-1982,}\/  
Progress mathematics {\bf 37}, Birkh\"auser Boston (1983) 401-426.
\vskip 0.2 cm

\item{[Miwa]}
T. Miwa, Painlev\'e property of monodromy preserving deformation 
equations and the analyticity of $\tau$-function, {\it Publ. RIMS,}\/ Kyoto 
Univ. {\bf 17} (1981) 703-721.
\vskip 0.2 cm

\item{[Nor]}
N.E. N\"orlund, The logarithmic solutions of the hypergeometric 
equation, {\it Mat. Fys. Skr. Dan. Vid. Selsk.}\/ {\bf 2}, no.5 (1963).
\vskip 0.2 cm

\item{[Ok]}
K. Okamoto,  Studies on the Painlev\'e equations I, sixth Painlev\'e 
equation, {\it Ann. Mat. Pura Appl.}\/ {\bf 146} (1987) 337-381.
\vskip 0.2 cm

\item{[Pain]}
P. Painlev\'e, Sur les Equations Differentielles du Second Ordre et
d'Ordre Superieur, dont l'Interable Generale est Uniforme, {\it Acta Math.}\/ 
{\bf 25} (1902) 1-86.
\vskip 0.2 cm

\item{[Pic]}
E. Picard, M\'emoire sur la th\'eorie des functions alg\'ebriques de deux 
varables, {\it Journal de Liouville}\/ {\bf 5} (1889), 135-319.
\vskip 0.2 cm

\item{[Sch]}
L. Schlesinger, \"Uber eine Klasse von Differentsial system beliebliger 
Ordnung mit festen kritischer Punkten, {\it J. fur Math.} {\bf 141} (1912), 
96-145.
\vskip 0.2 cm

\item{[Sib]}
Y. Sibuya, Linear differential equations in the complex domain: problems of 
analytic continuation, {\it AMS TMM}\/ {\bf 82} (1990).
\vskip 0.2 cm

\item{[SG]}
G. Sansone and J. Gerretsen, {\it Lectures on the theory of functions of a 
complex variable,}\/  P. Noordhoff editor, Groningen (1960).
\vskip 0.2 cm

\item{[Tak]}
L.A. Takhtajan, {\it Modular forms as tau-functions for certain integrable
reductions of the Yang-Mills equations,}\/ PREPRINT (1992)
\vskip 0.2 cm

\bye